\documentclass[poms,final,nonblindrev]{poms1_V1} 

\OneAndAHalfSpacedXI 

\usepackage{graphicx}
\usepackage{subfigure, epsfig}
\usepackage{natbib}

\usepackage{adjustbox} 
 
\DeclareMathOperator*{\argmaxB}{argmax} 
\usepackage{multirow}
\usepackage{rotating} 
\addcontentsline{toc}{section}{Appendices}

\newcommand{\vahid}[1]{{\color{red} #1}}
\newcommand{\mybox}{%
	\collectbox{%
		\setlength{\fboxsep}{1pt}%
		\fbox{\BOXCONTENT}%
	}%
}
\usepackage{ctable}
\usepackage{float}
\usepackage{subfigure}
\usepackage{makecell}
\usepackage[font=small,labelfont=bf]{caption}

 \bibpunct[, ]{(}{)}{,}{a}{}{,}%
 \def\bibfont{\small}%

\TheoremsNumberedThrough     
\ECRepeatTheorems

\EquationsNumberedThrough   
\begin{document}

\RUNAUTHOR{Doe and Doe}

\RUNTITLE{POM Template}

\TITLE{The Flood Mitigation Problem in a Road Network}

\ARTICLEAUTHORS{%
\AUTHOR{Vahid Eghbal Akhlaghi}
\AFF{H. Milton Stewart School of Industrial and Systems Engineering, Georgia Institute of Technology, Atlanta, GA, USA, \EMAIL{vahid.eghbal@gatech.edu}} 
\AUTHOR{Ann Melissa Campbell}
\AFF{Department of Business Analytics, Tippie College of Business, University of Iowa, Iowa City, IA 52242, United States \EMAIL{ann-campbell@uiowa.edu}}
\AUTHOR{Ibrahim Demir}
\AFF{Department of Civil and Environmental Engineering, University of Iowa, Iowa City, IA 52246, United States \EMAIL{ibrahim-demir@uiowa.edu}}
} 

\ABSTRACT{%
			Natural disasters are highly complex and unpredictable. However, long-term planning and preparedness activities can  help to mitigate the consequences and reduce the damage. For example, in cities with a high risk of flooding, appropriate roadway mitigation can help reduce the impact of floods or high waters on transportation systems. Such communities could benefit from a comprehensive assessment of mitigation on road networks and identification of the best subset of roads to mitigate. In this study, we address a pre-disaster planning problem that seeks to strengthen a road network against flooding. We develop a network design problem that maximizes the improvement in accessibility and travel times between population centers and healthcare facilities subject to a given budget. We provide techniques for reducing the problem size to help make the problem tractable.  We use cities in the state of Iowa in our computational experiments.
}%


\maketitle

\section{Introduction}

	Floods are more likely to occur than many other types of disasters \citep{sohn2006evaluating}. 
When a flood occurs, access to healthcare facilities is a top priority, followed by access to emergency operations facilities, such as supply distribution centers \citep{FEMA3}. When floods deteriorate and block transport systems, the loss of access to healthcare facilities (HCFs) creates challenges in providing continuing healthcare,  implementing a rapid vaccination, and transporting patients or injured people requiring treatment to hospitals \citep{mera2019towards,who,who2}. Thus, it is important to make decisions before such disasters occur regarding how to mitigate transportation infrastructure to preserve access to HCFs.

Investment in transportation networks plays an essential role in mitigation planning and decisions for many communities \citep{alabbad2022flood,teague2021collaborative}. Investments can strengthen the arcs of a network to enhance their survivability and the network's performance \citep{peeta2010pre}. 
Historically, local authorities often prioritized the selection of roads for mitigation based on either a simple heuristic approach or the road's usage, i.e., the volume of traffic \citep{bagloee2017identifying,mera2019towards}. To the best of our knowledge, in most of these approaches, roads are often selected without considering the impact they have on the rest of the network.

There are many approaches to mitigate the impact of flooding on communities, including elevating structures, relocation, and temporary flood walls \citep{yildirim2021integrated}.
In this work, we examine which roads and bridges to elevate to a level higher than the expected water level in large flood events.   We will refer to this practice as \textit{upgrading}.
To maintain the accessibility of the network's residents to healthcare centers and minimize their transportation time, it is essential to identify the roads with a risk of potential flooding, called \textit{vulnerable} roads, then upgrade the most critical ones to reduce the impact of flooding events \citep{alabbad2021assessment}. We assess the vulnerability of network segments in reference to the scenario of a 100-year flood. In an ideal scenario, all roads with the risk of potential flooding, called \textit{vulnerable} roads, would be upgraded. However, due to budget restrictions, not all of them can be mitigated at once.
This raises the question---how can city authorities prevent residents' loss of access to HCFs, and minimize the corresponding travel times by preparing a mitigation plan given budget restrictions?

To answer this question,
we define a transportation network that overlays a physical road network and enables the flow of patients and flood victims from population centers to HCFs. Due to the budget and facilities' capacity restrictions, we focus on selecting a subset of roads for mitigation to minimize the network's total travel time, such that each population center is served by a facility. We refer to this problem as the road network flood mitigation problem (RNFMP).
We develop a mathematical model to formulate RNFMP as an integer programming model with the aim of optimally allocating the mitigation budget to the roads and finding the optimal paths. We also identify the complexity of this new problem.

To the best of our knowledge, in all of the related work, the OD pairs (assignments of origins to destinations) are pre-specified or known in advance. This assumption separates our problem from the literature. Namely, the assignment of the origins to destinations must be determined, in addition to the selected arcs for mitigation in the RNFMP. Furthermore, for networks with pre-specified OD pairs, there exist lots of pruning techniques in the literature developed to reduce the size of the solution space and the number of decision variables, whereas
they are often inapplicable to RNFMP. 
Thus, our methodology involves both adapting existing pruning techniques and developing new ones to enable solving realistic problem instances. In this paper, we perform experiments with two cities in the state of Iowa with different-sized road networks.

The main contributions of this study are (1) providing a list of related topics (terms/keywords) in our literature review as a guiding framework for future studies; (2) formulating RNFMP without pre-specified OD pairs for the first time in pre-disaster arc-upgrading problems; (3) verifying the complexity of RNFMP; (4) introducing novel pruning techniques customized for RNFMP to solve large networks; and (5) offering a set of computational experiments to verify our approach for two cities in Iowa and a detailed case study of how the results can be used by road transportation authorities (e.g., DOT) in their decisions regarding a road network mitigation plan for one of the cities.

This paper is organized as follows: in Section \ref{lit}, a comprehensive literature review is provided. Section \ref{problem} provides the problem definition. Section \ref{imp3} introduces the developed model improvements, including the pruning techniques and valid inequalities.  Section \ref{DOE} introduces our dataset, representing two road transportation networks of two cities located in the State of Iowa, and describes its features for experimental design. Section \ref{impimp} utilizes this dataset to evaluate the effectiveness of our proposed improvements developed in Section \ref{imp3}, and Section \ref{results}  conducts a sensitivity analysis for RNFMP's parameters and reports the results for our case study.  Finally, we provide concluding remarks in Section \ref{con} and discuss future work.

\section{Literature Review}\label{lit}
Due to the massive impact that natural disasters can have on the health of the residents of an entire region \citep{alabbad2022comprehensive}, the study of ways to better understand and mitigate such incidents has become an area of rising concern \citep{bagloee2017identifying}. Specifically, the failure of road segments in transportation networks has been an important research topic in transportation planning since the 1990s \citep{santos2010interurban}. A growing interest has been observed in research related to \textit{resilience} and \textit{vulnerability} 
of transportation networks \citep{nagurney2012fragile} because the selection of a subset of critical roads for mitigation has usually been determined within the context of these two subjects/terms. However, the literature has not provided consistent definitions for them \citep{faturechi2014travel,woods2015four}.  

\cite{mattsson2015vulnerability} provide a review of research on the \textit{vulnerability} and \textit{resilience} of transportation networks and show that they have been defined together by related concepts.  \cite{berdica2002introduction} defines vulnerability as susceptibility to catastrophes that can cause reductions in transportation network serviceability, whereas \cite{bruneau2003framework} define resilience as the system's ability to absorb a shock or an abrupt reduction of performance after a shock. From here on, we adopt these definitions when we speak of a network's vulnerability or resilience. Inspired by these definitions, we measure the performance of our network based on the change in the network's total transportation time after a flood event. For an extended review of definitions and measures, we refer the reader to  \cite{bruneau2003framework} and \cite{yucel2018improving}.

In Section 2.1, we provide a detailed review of the most relevant studies and discuss where our study fits in the literature.
Finally in Section 2.2, we briefly introduce the related terms/topics corresponding to the problems with similar settings to RNFMP.

\subsection{The positioning of this paper in the relevant literature}\label{rel-papers}
The related literature can be categorized into two groups based on two distinct streams: (i) papers that are rooted in graph theory and study the vulnerability of networks based on their topological (i.e., the spatial layout of the network) properties, (ii) papers that represent the demand and supply side of the transportation network to allow for a complete assessment of the consequences of disruptions for the users and society \citep{bagloee2017identifying}. While the first group is only concerned with the topological characteristics of networks, such as shortest paths \citep{mount2019towards}, centrality, transitivity, etc., the latter considers the dynamic characteristics of flow in networks (e.g., the flow of flood victims', HCFs' capacities, and the number of residents (weight) associated with each population center studied in the current paper). In the first group, the proper functionality of the system and the criticality levels of roads have been evaluated based on the network's connectivity \citep{ulusan2018restoration,bagloee2017identifying}.
However, in the second group,
a highly fragmented network that is split into distinct components can still be fully functional such that the flows are unaffected \citep{morris2012transport}. 
For a recent review of these traditions, see \cite{bagloee2017identifying} and \cite{ulusan2018restoration}.

As already highlighted, the decisions about the selection of the set of critical roads should consider its impact on mitigating the increase in network's transportation time. Also, according to \cite{faturechi2014travel} and \cite{bagloee2017identifying}, when the flow dynamics are important as the overall network performance, the total travel time is a widely used index. As a result, in the current study, the total travel time has been used as RNFMP's objective function.

We provide a detailed review of some of the most relevant and recent studies in the literature.
The main differentiating point of our work with these papers is that our model finds the optimal mitigation investments and optimal OD assignments simultaneously.  But in the related papers, OD assignments are pre-specified as the model's parameters where each origin is paired with a specific destination, usually the nearest one. 
Subsequently, they often don't consider the corresponding demands or capacity restrictions.
Unlike the literature, our model does \textit{not} limit each origin to be served by a pre-specified destination. One reason for adding this complicating assumption is the role of cellphones and social media as effective tools for local, state, and federal governments to communicate with the public during catastrophes and share 
quick updates about road closures and functional routes to hospitals or shelters \citep{yan2019social,FEMA2,sit2019identifying}. 

In February 2022, FEMA released an updated texting feature, available across all U.S. states and territories, that allows public to find shelter addresses during a disaster. ``Disasters frequently disrupt communications systems which can leave survivors feeling overwhelmed and helpless when they are trying to locate shelters," said FEMA Administrator Deanne Criswell. ``Since texting capabilities are often unaffected during disasters, our updated Text to Shelter option is an easy and accessible way survivors can locate nearby shelters with a tap of a button." Through their collaboration with the American Red Cross, users can click on shelter addresses provided inside the text message and view the directions on their phones \citep{femaredcross}.

Effective information exchange is essential for a successful humanitarian response as it enables the provision of relevant and timely information to support appropriate decision-making \citep{altay2014information}. \cite{yan2019social} study the impact of this information exchange on social engagement during disaster preparedness, response, and recovery using Facebook data from five benchmark organizations that responded to Hurricane Sandy in 2012. In the past, when officials issued paper news releases, it was sensible to assume that residents would drive to their \textit{closest} available HCF in their community. However, in today's digital age, social media platforms have become a key tool for communication during natural disasters.	According to recent studies, more than 85 percent of local governments use social media platforms to communicate with their constituents \citep{FEMA2}.

Since social media has increasingly become a key tool for communication during natural disasters, community officials can assign each residential area to a specific HCF and guide them with post-disaster functional routes and directions. The 2015 Historic South Carolina Flood is a great example highlighting how social media can help critical information (e.g., transportation information) flow much faster than it did in the past during an emergency \citep{carolina}. To the best of our knowledge, our study is the first that considers both the OD assignments and road upgrading investments as model's decision variables.
Below, we review the most related and recent studies.

\cite{peeta2010pre} consider a pre-disaster planning problem to strengthen a highway network whose arcs are subject to random failures due to a disaster. They assume that each road's failure probability is known, and mitigation investment decreases the likelihood of failure. Unlike our paper that considers real mitigation costs for roads and bridges, they consider the number of roads, where each arc's mitigation cost is equal to one unit. In Section \ref{results}, we will show the drawbacks of this approach. They consider the network connectivity for first responders between pre-specified OD pairs. The goal is to select the best set of arcs to strengthen under a limited budget to maximize the post-disaster connectivity and minimizing traversal costs between the OD pairs. The problem is modeled as a two-stage stochastic program. Since the OD pairs are known, first, they minimize the expected cost of OD paths by choosing the arcs to be upgraded. If an OD path is disconnected, a high value of (penalty) cost will be considered for it. Second, they minimize the total travel time. The computational study is based on the highway networks of Istanbul, with 20 nodes and 30 arcs.

Accessibility was first defined by \cite{hansen1959accessibility} as the number of potential opportunities (destinations) that can be reached by each location (origin) within a given travel distance or travel time. By adopting this definition, \cite{santos2010interurban} focus on a road network design problem to determine the best investment decisions to be made regarding road improvements. They propose a bi-objective optimization model to maximize the accessibility and robustness objectives simultaneously. To address the robustness criterion, they apply three measures through spare network capacity, evacuation capacity, and vulnerability. Commuters are assumed to follow the routes that minimize their costs. They develop a genetic algorithm to solve large instances and discuss the results obtained through the model for three randomly designed networks.

\cite{tong2015transportation} aim to increase users' accessibility to desired destinations or between major activity locations, subject to a given construction (not mitigation) budget and a given travel time budget. As a result, the objective is to minimize the inaccessibility among all major activity locations, where a one-unit inaccessibility cost ($c_{ij}=1$) will be added to the objective if location \textit{j} cannot be reached from location \textit{i} on time. They utilize a Lagrangian relaxation method to solve their proposed MILP model. The computational performance of their developed solution approach has been evaluated on the large-scale Chicago sketch network with 933 nodes, 2,950 arcs (including 149,382 given OD pairs), and a set of 20 candidate arcs to be built subject to a construction budget of five arcs.

In another related paper, \cite{angulo2016lagrangian} address expanding an existing highway network by creating new roads to minimize total travel time. 
With a given budget, instead of finding the optimal set of roads to be upgraded, their goal is to construct new roads between cities that are not directly connected. The construction cost of each road is equal to one unit, and similar to the rest of the literature, OD pairs are pre-specified.  
The case study is carried out for the Castilla-La Mancha (Spain) region with 2,154 roads and 290 nodes. They solve their model using CPLEX, where the stopping criterion was a gap of less than 4\%. 

 \cite{bagloee2017identifying} identify the most critical roads to mitigate and increase their resilience to minimize total travel time. First, they apply the sensor (loop detector) location problem (SLP), within which a selected number of high-demand roads are heuristically identified as ``candidate" critical roads.  Unlike our paper that finds the optimal set of roads to be mitigated, they find the set of critical roads by calculating each road's demand (traffic) based on the given OD demands. Then they find the critical combination scenarios by solving a series of discrete network design problems (DNDP). The DNDPs are solved based on an optimal relaxation method using Bender's Decomposition with a 2\% gap. 
Unlike our paper that calculates the real mitigation cost values in dollars, the mitigation cost of each arc is equal to one unit. They undertake the network of the city of Winnipeg (Canada) as a case study with 943 nodes and 3075 arcs.

\cite{yucel2018improving} address a network improvement problem to improve the network's resiliency against disasters. The main goal is to improve the \textit{expected} post-disaster functionality of a highway network by optimizing the pre-disaster road improvements (strengthening decisions). First, to predict the post-disaster status of the network, they apply a dependency model for random arc failures. Then, they estimate the accessibility measure by checking pre-generated short and dissimilar paths in the sample. This approach necessitates a two-stage stochastic programming framework.  A tabu search (TS) algorithm is used to solve the problem. The computational analysis is applied to a case study of Istanbul under the risk of an earthquake to obtain insights for preparedness activities. The network consists of 60 nodes and 83 undirected arcs.

 \cite{mera2019towards} address the problem of road maintenance and development to improve a network's resilience. The available capacity of a road is prone to be affected by a natural disaster leading to a drop in its capacity. Standard road capacities are given, and disrupted arcs have an available capacity between 0 and 50\% of the standard capacity, whereas upgraded arcs can have a capacity up to an arbitrary value of 120\%. Consequently, the serviceability of an arc reduces when an incident adversely affects the capacity of the arc. The serviceability index is defined as a fraction of the total available capacity to standard capacity. This definition indicates that while multi-lane (higher-order) road arcs are key to maintaining the network functionality, single-lane (lower-order) road arcs are equally critical in contributing to the network connectivity and could play a significant role in reducing the network's disruption. However, in their experiments, the authors assume that all of the streets have two lanes. The objective is to minimize a measure of vulnerability to disruption under budgetary constraints. A Simulated Annealing Metaheuristic approach is used to solve the problem. To test its performance, an instance based on a segment of the City of York in England is generated with 110 nodes and 251 arcs.

There are several studies that address the network restoration and recovery problems in the post-disaster stage to maximize various measures, such as accessibility.  For a comprehensive review of post-disaster network improvement problems, we refer the reader to \cite{ccelik2016network}.

\subsection{Related terms/problems with similar settings}
There exist several related terms and topics, such as \textit{network strengthening problem (NSP)}, \textit{best upgrade plan problem (BUP)}, \textit{network upgrading problem (NUP)}, etc., which correspond to problems with settings similar to RNFMP. 
In this section, we briefly review these terms/keywords and show that, although they are called by different names, they mostly address the same problem. However, we believe that outlining these terms provides a guiding framework for future studies. 

\cite{peeta2010pre} propose a two-stage stochastic program to increase the post-disaster survival probability of arcs by strengthening them in the pre-disaster stage. The objective of this problem, which was later called the \textit{network strengthening problem (NSP)} by  \cite{yucel2018improving}, is to minimize the shortest paths between pre-specified OD pairs. The setting of this problem is similar to that of RNFMP, but minimizing the total travel time in RNFMP is not identical to minimizing the shortest paths in NSP, where (i) OD pairs are known; (2) there are no capacity restrictions; and (3)  travel times are not weighted based on the populations.

The resource constrained \textit{best upgrade plan problem (BUP)} that has been studied by  \cite{lin2015best} selects a subset of arcs among the upgradable ones so that the shortest paths between ODs are minimized, and at the same time, the summed upgrade cost does not exceed a specific budget. There exists a cost associated with each arc that, if spent, the weight (e.g., traveling time) of the arc can be reduced to a new value. \cite{lin2015best} highlight that the closest problem to BUP is the \textit{network upgrading problem (NUP)}. In NUP, the aim is to determine which arcs should be upgraded to which levels to improve the network's performance, subject to budget restrictions \citep{yucel2018improving}.  

\cite{duque2013accessibility} study the accessibility measure by presenting the \textit{accessibility arc upgrading problem (AAUP) }to determine which arcs should be upgraded under budget restrictions such that the total weighted travel time from a set of demand nodes to their \textit{closest} supply centers is minimized. They utilize the accessibility concept such that the accessibility of a node is measured by its travel time to its \textit{closest} center, which makes the set of OD pairs pre-specified. Two heuristic approaches are proposed to deal with this problem. The authors show that AAUP is a variant of the minimum cost flow problem containing additional decision variables corresponding to the upgrading decisions. Due to budget restrictions, these problems are also referred to as \textit{budget constrained network upgrading problems (BC-NUP)} \citep{krumke1999improving}.

\cite{campbell2006upgrading} define the \textit{\textit{q}-upgrading problem }that  identifies \textit{q} arcs in the network to upgrade such that the travel time of the maximum shortest path is minimized. upgrading an arc connecting two nodes corresponds to using a faster mode of transportation, such as a plane instead of a truck. As a result, similar to BUP and NUP, when an arc is upgraded, its travel time decreases. This problem, which is NP-Complete, determines which arcs should be upgraded such that the maximum shortest path between any OD pair in the network is minimized under budget restrictions. In arc upgrading models, three cases are considered depending on the values to which each arc can be upgraded: binary (upgraded, non-upgraded), integer (several upgrading levels), rational (continuous values in an interval). For the binary upgrades, the authors show that minimization of the total weight of the graph under the budget constraint is NP-hard and propose three heuristic algorithms to solve this problem.

\cite{demgensky2002flow} find the optimal arcs to upgrade under budget restrictions for a minimum cost flow problem. They show that the problem is NP-hard and propose an approximation algorithm to solve it. For more examples and applications of UAPs, we refer the reader to \citep{duque2013accessibility,lin2015best,nepal2009upgrading}. In Section \ref{NP}, we will determine the complexity of the RNFMP in detail.

The general \textit{network design problem (NDP)} is to find the best choice of affordable candidate projects while accounting for the way commuters utilize the network \citep{bagloee2017identifying}. \cite{tong2015transportation} define the \textit{transportation network design problem (TNDP)}, where the aim is to improve the transportation network's performance through introducing new arcs or improving existing ones. The objective is to minimize inaccessibility in terms of unreachable activity locations within given travel time and construction budget constraints. \cite{gutierrez1996robustness} propose the \textit{uncapacitated network design problem (U-NDP)} to determine the best configuration of the network that minimizes the sum of the fixed costs of arcs chosen to be in the network along with the costs of routing goods through the network defined by these arcs. 

\cite{santos2010interurban} categorize NDPs into two groups: \textit{discrete network design problem (DNDP)} and \textit{continuous network design problem (CNDP)}. The former explores the addition of new arcs or binary capacity improvement decisions in a road network, which is similar to the binary arc upgrading models in NUP. On the other hand, CNDP focuses on the continuous expansion of the capacity of existing arcs \citep{mera2019towards}. We refer the interested reader to \citep{farahani2013review} for more examples and applications.

\cite{chen2002capacity} claim that when the disruption is a huge disaster rather than simple congestion, the binary arc upgrading approach is most rational. Even though there are studies with continuous or discrete approaches, where a complete arc shutdown may not occur if the inundation is not severe,
reports presenting flood-related incidents indicate that fatalities are quite commonly caused by vehicles entering road sections, giving a false impression of not being overly flooded \citep{wisniewski2020vulnerability}, given that it only takes 12 inches of water to carry away a small car \citep{FEMA2}. Since we assume the operational routes to the HCFs will be given to commuters by local governments, it makes sense that community officials don't want to route them through roads that may look  partially functional with possible false impressions. Thus, the binary approach has been selected for RNFMP.

\section{Problem Characteristics}\label{problem}
We aim to strengthen a road network in a given city (geographical area) in which a set of population centers along with a set of HCFs is considered. Section \ref{def} introduces the sets, parameters, and variables for our problem. Section \ref{model3} proposes our formulation, and Section \ref{NP} proves the NP-hardness of the RNFMP.
\subsection{Problem definition and notation}\label{def}
The existing network is modeled by a directed graph $\mathbb{G}(\mathcal{N},\mathcal{A})$ in which $N=\{\mathcal{N}_o\cup \mathcal{N}_\tau \cup \mathcal{N}_d\}$ is a node-set and $\mathcal{A}=\{\mathcal{A}_v,\mathcal{A}_n\}\subseteq\mathcal{N}\times\mathcal{N}$ is an arc set. The nodes in $\mathcal{N}_o$ are called origins (population centers), while the nodes in $\mathcal{N}_d$ are called destinations (HCFs), and $\mathcal{N}_\tau$ is the set of transshipment nodes, where $\mathcal{N}_o\cap \mathcal{N}_\tau \cap \mathcal{N}_d\in \emptyset$. Arcs in $\mathcal{A}$ represent roads, and the subset of arcs $\mathcal{A}_v$ is composed of all vulnerable roads that will be non-functional (flooded) after a flood event and can be upgraded. However, non-vulnerable arcs in $\mathcal{A}_n$ will be operative, where $\mathcal{A}_v \cap \mathcal{A}_n \in \emptyset$. 

Road segments' vulnerabilities are determined based on the scenario of a 100-year flood map, a commonly applied standard in flood risk assessment  \citep{ewing2021ethical,sermet2020serious,gutry2006numerical}.
The 100-year flood maps utilized in this study are obtained from the Iowa Flood Center.
In a 100-year flood map, if any road is within a 100-year floodplain (i.e., within the areas of low-lying ground subject to flooding), it is assumed to have a potential risk of failure when a flood occurs \citep{fema2016definitions} and thus will be considered in our network's set of vulnerable roads. Complete shutdowns have been considered for those road segments damaged by the flood, which is also the most common approach in the literature \citep{wisniewski2020vulnerability}. Following \cite{yucel2018improving}, it is assumed that only the arc segments are subject to failure, but not the nodes, as the risk associated with nodes are incorporated into the arcs incident to them. 

The directed arc $(i,j)$ will be included in the set of arcs, i.e., $(i,j) \in \mathcal{A}$, if there exists a road that connects nodes $i$ and $j$ without going through any other node. For each arc $(i,j)\in \mathcal{A}$, let $t_{ij}\ge0$ denote travel time, and for each arc $(i,j)\in \mathcal{A}_v$, let $c_{ij}\ge0$ denote its mitigation cost. The travel time $t_{ij}$ is a function of $(i,j)$'s length and speed limit. There is a total mitigation budget $B$. For each node $k$ in $\mathcal{N}_o$, a weight $w^k$ is defined as a function of the number of residents associated with that node ($h^k$) and corresponds to its weight in the objective function. This gives a higher priority to the more populated areas. Each node $j$ in $\mathcal{N}_d$ is associated with a capacity $H_j$, taking into account the destination's capacity. 

Some types of network flow problems split the flow from a particular origin.
In this study, every resident at the same origin is assigned to the same destination. Each origin represents a group of people located in one small area. Community officials send quick updates to the residents about road closures and functional routes, guiding them to the specific medical centers they should go to. A mitigation plan/solution specifies a set of arcs among $\mathcal{A}_v$ to be upgraded such that the weighted sum of the travel times required to travel from each origin $i$ in $\mathcal{N}_o$ to a destination $j$ in $\mathcal{N}_d$ in the upgraded network is minimized. A mitigation solution is feasible if the total cost of upgraded roads does not exceed budget $B$, and each origin is served by a destination. Namely, feasible solutions to the RNFMP consist of: (i) a subset of arcs to be upgraded represented by $y_{ij}$'s; (ii) a subset of arcs that connect each origin to its allocated destination represented by $x_{ij}^k$. The sets, parameters, and decision variables used to formulate RNFMP are summarized in Table \ref{notation}.
\begin{table}[!ht]
	\tiny
	\caption{Notation used for RNFMP}
	\centering
	{\def\arraystretch{1} 
\begin{adjustbox}{max width=\textwidth}				
		\begin{tabular}{rll}
			\hline
			\hline
		\multicolumn{2}{l}{\textbf{Sets and indices}} &  \\
& $\mathcal{A}$: &  The set of all arcs \\
& $\mathcal{N}$: &  The set of all nodes that serve as endpoints of arcs in $\mathcal{A}$, indexed by $i$ or $j$ \\
& $\mathcal{N}_o$: &  The set of origins (population centers), indexed by $k$ \\
& $\mathcal{N}_d$: &  The set of destinations (HCFs) \\
& $\mathcal{N}_\tau$: &  The set of transshipment nodes \\
& $\mathcal{A}_n$: &  The set of non-vulnerable arcs \\
& $\mathcal{A}_v$: &  The set of vulnerable arcs \\
&       &  \\
\multicolumn{2}{l}{\textbf{Parameters}} &  \\
& $t_{ij}$: & The travel time of arc $(i,j)\in \mathcal{A}$ \\
& $c_{ij}$: & The cost of mitigating arc $(i,j)\in \mathcal{A}_v$ \\
& $B$     & The total available budget for mitigation \\
& $h^k$: & The population of origin $k \in \mathcal{N}_o$ \\
& $H_j$: & The capacity of destination $j \in \mathcal{N}_d$ \\
& $w^k$: & The weight of origin $k \in \mathcal{N}_o$ that is a function of $h_k$ \\		
&       &  \\
\multicolumn{2}{l}{\textbf{Decision Variables}} &  \\
& $x_{ij}^k$: & Equal to 1 if arc $(i,j) \in \mathcal{A}$ is traversed by residents associated with origin $k$; 0 otherwise \\
& $y_{ij}$: & Equal to 1 if arc $(i,j)\in \mathcal{A}_v$ is mitigated; 0 otherwise \\
\hline\hline
		\end{tabular}
\end{adjustbox}	
		\label{notation}
	}
\end{table}
\subsection{The mathematical formulation of RNFMP}\label{model3}
Below, a mathematical integer programming formulation is presented to model RNFMP:
{\footnotesize
		\small
\begin{align}
& \mbox{Min }  \sum_{k \in \mathcal{N}_o} \sum_{(i,j) \in \mathcal{A}} t_{ij}w^k x_{ij}^k\;  \label{1}\\
\nonumber & \mbox{Subject to}\\
&\sum_{i:(i,k) \in \mathcal{A}}x_{ik}^k - \sum_{h:(k,h) \in \mathcal{A}}x_{kh}^k =-1 & \forall k \in \mathcal{N}_o \label{2}\\
&\sum_{i:(i,j) \in \mathcal{A}}x_{ij}^k - \sum_{h:(j,h) \in \mathcal{A}}x_{jh}^k \ge 0 & \forall k \in \mathcal{N}_o, \forall j \in \mathcal{N}_d \label{3}\\
&\sum_{i:(i,j) \in \mathcal{A}}x_{ij}^k - \sum_{h:(j,h) \in \mathcal{A}}x_{jh}^k =0 & \forall k \in \mathcal{N}_o, \forall j \in \mathcal{N}\setminus\mathcal{N}_d, k\ne j \label{4}\\
& \sum_{(i,j) \in \mathcal{A}_v} y_{ij}*c_{ij} \le B \label{5}\\
& x_{ij}^k\le y_{ij} &\forall k \in \mathcal{N}_o, \forall (i,j) \in \mathcal{A}_v \label{6}\\
&y_{ij} \le \sum_{k \in \mathcal{N}_o} x_{ij}^k  &\forall (i,j) \in \mathcal{A}_v \label{7}\\
&\sum_{k \in \mathcal{N}_o} \sum_{i:(i,j) \in \mathcal{A}}x_{ij}^k h^k - \sum_{k \in \mathcal{N}_o} \sum_{h:(j,h) \in \mathcal{A}}x_{jh}^k h^k\le H_{j}\; &\forall j \in \mathcal{N}_d \ \label{8}\\
& x_{ij}^k \in \{0,1\}\; &\forall k \in \mathcal{N}_o, \forall (i,j) \in \mathcal{A} \label{9}\\
& y_{ij} \in \{0,1\}\; &\forall (i,j) \in \mathcal{A}_v \label{10}
\end{align}
}

The objective function \eqref{1} minimizes the sum of the weighted travel times.  Constraints \eqref{2}-\eqref{4} are the usual network flow conservation constraints. Note that each destination can serve more than one origin. Constraints \eqref{5} are typical knapsack constraints, which restrict the total mitigation investment to the available budget $B$. Constraints \eqref{6} limit a flooded road to be used only if it is mitigated. Constraints \eqref{7} prevent a flooded road from being mitigated if it is not being used. Constraints \eqref{8}, which are known as the supply nodes' mass balance constraints in the \textit{min-cost flow problem (MCFP)}, ensure that the total population served at a destination does not exceed its capacity. Constraints \eqref{9} and \eqref{10} restrict $x_{ij}^k$ and $y_{ij}$ to zero or one.

\subsection{Complexity}\label{NP}
To prove that RNFMP is NP-hard, we will reduce the decision version of the \textit{generalized assignment problem (GAP)}, which is proved to be NP-hard \citep{fisher1986multiplier} to the RNFMP's decision problem. 
Let $\mathbb{G}(\mathcal{N},\mathcal{A})$ be a directed graph where $\mathcal{N}$ and $\mathcal{A}$ are defined as indicated in Table \ref{notation}. 
Let's define $Z$ as the sum of the weighted travel times, in which the travel time of each arc $t_{ij}$ is multiplied by a weight $w^k$ representing the number of residents associated with origin $k$. 

To prove the NP-hardness of RNFMP, we consider a case where the mitigation budget is large enough to upgrade all vulnerable arcs (i.e., $B\ge \sum_{(i,j) \in \mathcal{A}_v} c_{ij}$). This special case represents a simpler version of RNFMP, in which $y_{ij}$ variables can be set to 1 for every $(i,j)\in \mathcal{A}_v$. Then, the question in the decision version of RNFMP is: Given $\mathbb{G}$, with $h^k\ge 0$ $\forall k\in \mathcal{N}_o$, $H_j\ge 0$ $\forall j \in \mathcal{N}_d$, $t_{ij}\ge 0$ $\forall (i,j)\in \mathcal{A}$, and a positive rational $\beta$, does there exist a weighted sum of travel times $Z$ that has a value less than or equal to $\beta$?

Note that when $y_{ij}$'s are set to 1, the travel time from an origin $k \in \mathcal{N}_o$ to a destination $j \in \mathcal{N}_d$ will be equal to the value of the shortest path from $k$ to $j$. By finding the shortest path for every OD pair, which can be solved by polynomial-time algorithms such as Dijkstra's algorithm, we can remove the set of arcs $\mathcal{A}$ from $\mathbb{G}$, and add an arc from every origin to every destination representing the shortest path for that specific OD pair. Subsequently, $x_{ij}^k$ variables can be replaced by $x_{j}^k$'s, where $x_{j}^k=1$, if origin $k$ is assigned to destination $j$, 0 otherwise. 

The GAP considers the minimum cost assignment of $n$ jobs to $m$ agents. Each agent is assigned to one or more jobs, requiring some resources, such that no agent exceeds its capacity. By defining $Z^\prime $ as the cost of assigning jobs to agents, and a positive rational $\beta^\prime$, the decision version of GAP is to verify if there exists a feasible jobs-agents assignment in which $Z^\prime \le \beta^\prime$?

By letting (i) $h^k$ that equals the population of each origin in RNFMP, be the number of resources required for job $k$ in GAP, (ii) $t_{ij}$, which represents the shortest path from $i$ to $j$ in RNFMP, be the cost of assigning job $i$ to agent $j$ in GAP, and (iii) $w^k$ be equal for all of the origins (e.g., $w^k=1 \forall k\in \mathcal{N}_o$), RNFMP can be formulated as a GAP problem:

{	\small
\begin{align}
\nonumber& \mbox{Min }  \sum_{k \in \mathcal{N}_o} \sum_{j \in \mathcal{N}_d} t_{kj} x_{j}^k\;  \\
\nonumber & \mbox{Subject to}\\
\nonumber&\sum_{k \in \mathcal{N}_o} x_{j}^k h^k \le H_{j}\; &\forall j \in \mathcal{N}_d \ \\
\nonumber&\sum_{j \in \mathcal{N}_d} x_{j}^k =1\; &\forall k \in \mathcal{N}_o \ \\
\nonumber& x_{j}^k \in \{0,1\}\; &\forall k \in \mathcal{N}_o, \forall j \in \mathcal{N}_d 
\end{align}}
Thus, RNFMP and GAP are polynomially equivalent. The RNFMP can be transformed into a GAP problem, and each instance of the GAP problem can be transformed into an instance of the RNFMP.

\section{Model Improvements}\label{imp3}
In Section \ref{model3}, we proved the RNFMP to be an NP-hard problem. Our initial computational experiments confirm that the RNFMP is difficult to solve for realistic instances. 
This section discusses several improvements used to help reduce the solution time without compromising optimality. In Section \ref{impimp}, we will demonstrate the computational impact of the ideas from this section.

In Section \ref{PT}, we develop several techniques to prune a network's nodes and arcs. Section \ref{VIVR} introduces several techniques that reduce the number of decision variables without pruning the network. Finally, in Section \ref{init}, we develop a greedy approach to determine the initial solution for RNFMP.
\subsection{Network pruning techniques}\label{PT}
We propose several techniques to prune $\mathbb{G}$. We use $\delta_{in}(i)$ and $\delta_{out}(i)$ notation to denote the set of node $i$'s incoming and outgoing arcs, respectively. We introduce $\mathcal{N}_{p}$ to be the set of articulation points in $\mathbb{G}$. An articulation point (AP) is a node whose removal disconnects and divides the network into $n\ge2$ components. For every node $v\in \mathcal{N}_p$, we let $\Delta(\mathcal{N}_{\Delta},\mathcal{A}_{\Delta})$ represent one of the connected components of the subgraph $\mathbb{G}-v$. We will leverage the APs, which can be identified in $O(|\mathcal{N}|+|\mathcal{A}|)$ linear time \citep{tarjan1972depth}, to potentially prune connected components.

In Technique 1, we use articulation points (AP) to prune potential nodes and arcs which belong to specific $\Delta$ components. A related AP-based idea was used by \cite{beldiceanu2005tree} for reducing the number of variables associated with tree constraints in graphs.

\begin{figure}[h]
\begin{center}
		\caption{Pruning techniques}
	\includegraphics[scale=0.2]{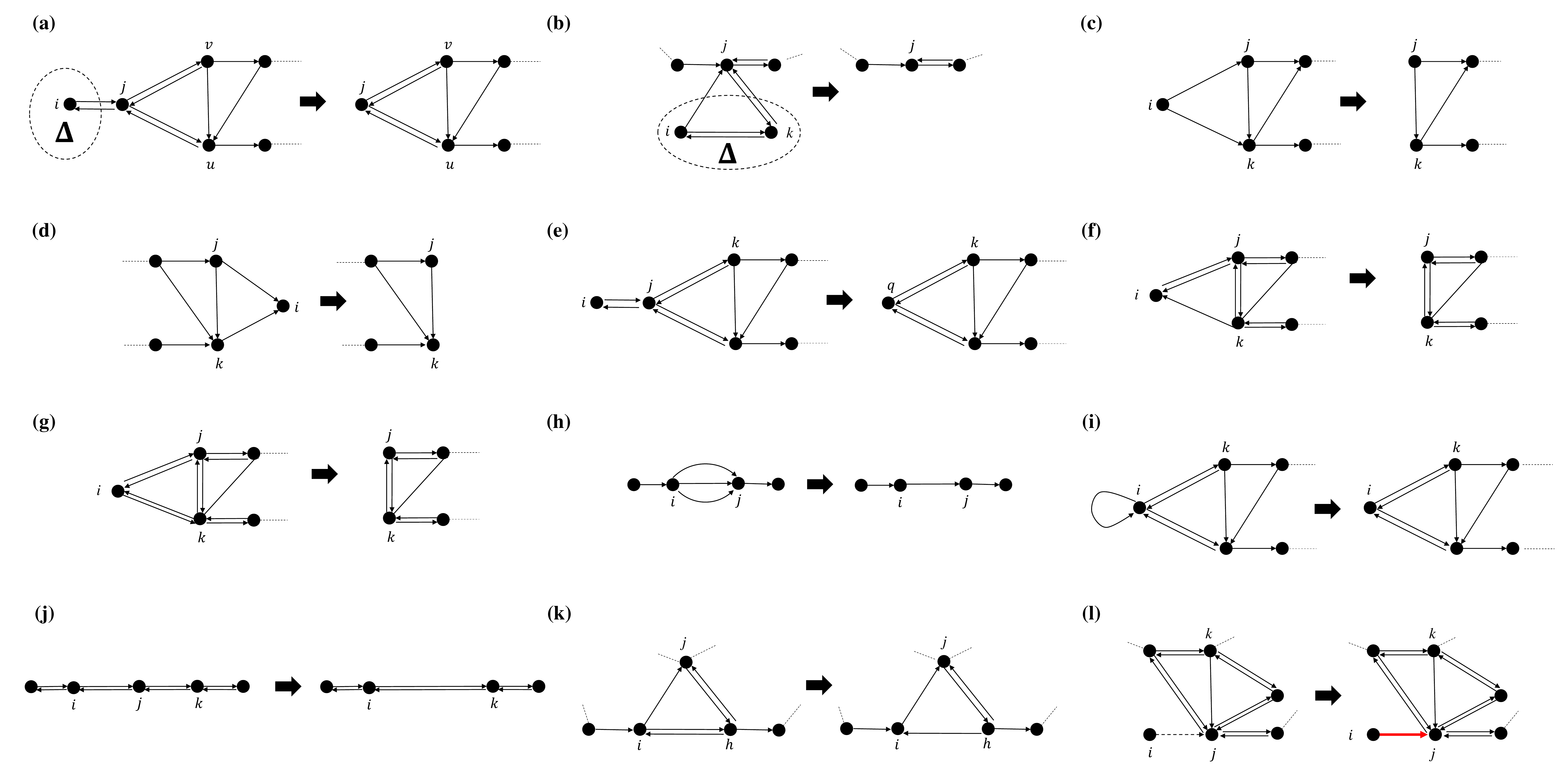}
	\label{tech}
\end{center}
\end{figure}

\begin{itemize}
	
	\item[]{\textbf{Technique 1:}\textit{
			For every node $j\in \mathcal{N}_{p}$, a disconnected component $\Delta$ obtained from $\mathbb{G}\setminus\{j\}$ can be pruned if $\nexists i\in \mathcal{N}_{\Delta} : i \in \mathcal{N}_o \cup \mathcal{N}_d$.
		}
		
		By definition of articulation points, any path between node $i \in \mathcal{N}_{\Delta}$ and node $v \notin \mathcal{N}_{\Delta}\cup \{j\}$ contains node $j$. Since $\Delta$ contains no origin (destination), there is no path initiating from (ending at) $\Delta$. Thus, any resident of the network, who enters component $\Delta$ through node $j$ (from an origin located outside $\Delta$) needs to leave $\Delta$ by revisiting node $j$ to reach a destination. This creates a circuit on this path with a start and end at node $j$, which is sub-optimal.  For instance, assuming there exists only one node $i$ in $\Delta$, path  $v\rightarrow \textbf{j} \rightarrow \textbf{i} \rightarrow \textbf{j }\rightarrow u$ is sub-optimal compared to $v\rightarrow \textbf{j} \rightarrow u$, where $v,u \notin \mathcal{N}_{\Delta}$ (Figure \ref{tech}a and \ref{tech}b). 
	}
\end{itemize}	
Techniques 2 through 4 are inspired by the tree pruning heuristic (TPH) technique developed by  \cite{bi2014graphical} to prune redundant children from trees. Though Techniques 2 through 4 are adapted and customized in three different ways given the specific characteristics of the RNFMP.
\begin{itemize}
	\item[]{\textbf{Technique 2:}\textit{
			Any node $i\in \mathcal{N}_\tau$ with either no incoming arc ($\delta_{in}(i)=0$) or no outgoing arc ($\delta_{out}(i)=0$) can be eliminated from $\mathbb{G}$.}
		
		If the transshipment node $i$ does not have either incoming or outgoing arcs, there can be no flow through node $i$ . As a result, node $i$ with its outgoing or ingoing arcs can be removed without affecting the solution (Figures \ref{tech}c and \ref{tech}d). 
		
	}

	\item[]{\textbf{Technique 3:}\textit{
			Adjacent nodes $i$ and $j$ can be merged if $i\in \mathcal{N}_o$, $j \in \mathcal{N}_\tau$, and there exist non-vulnerable arcs $(i,j)$ and $(j,i)$ where $\delta_{in}(i)=\delta_{out}(i)=1$.} 
		
		As mentioned in Section \ref{problem}, $\mathcal{N}_o$ and $\mathcal{N}_d$ are exclusive sets. Similar to Technique 1, since node $i$ does not belong to the set of destinations and
		is connected to the rest of the network only though node $j$, it is not optimal for the residents associated with any other origin $k\ne i$ to traverse node $i$. On the other hand, for origin $i$ to be served or assigned to a destination, its residents have no option other than crossing node $j$. Thus, it is possible to eliminate arcs $(i,j)$ and $(j,i)$, and merge nodes $i$ and $j$ into a single node $q$ (Figure \ref{tech}e), which will be added to set of origins. The scalar value of the corresponding weighted travel time from $i$ to $j$ can be added to the objective function in the pre-processing step of the solution approach.}

\item[]{\textbf{Technique 4:}\textit{
		If an isolated non-vulnerable triangle is established on three nodes $i\in \mathcal{N}_\tau$ and $j,k \in \mathcal{N}$, where $i$ is only connected to the rest of $\mathbb{G}$ through nodes $j$ and $k$, node $i$ and every arc connecting $i$ and $j$ or $k$ can be removed if $t_{ki}+t_{ij}\ge t_{kj}$  (Figures \ref{tech}f and \ref{tech}g) and $t_{ji}+t_{ik}\ge t_{jk}$  (Figure \ref{tech}g).}  
	
	The explanation is based on the properties of triangle inequality.
	
}

\end{itemize}
Technique 5 is related to the multiple-path pruning algorithm, proposed by \cite{poole2017python}. 
\begin{itemize}
\item[]{\textbf{Technique 5:}\textit{
		If parallel non-vulnerable arcs exist between nodes, only the arc with the shortest travel time can be kept while the rest can be eliminated, so that there is a single arc connecting each pair of nodes.} 
	
	If two or more arcs are incident on precisely the same nodes having the same direction, $\mathbb{G}$ can be pruned by eliminating the one(s) with a longer travel time (Figure \ref{tech}h). }
\end{itemize}
Technique 6 is inspired by the cycle pruning algorithm, proposed by \cite{poole2017python}. Given the specific setting of the RNFMP, we have customized \textit{cycle pruning} algorithm proposed to prune \textit{loops}.
\begin{itemize}
\item[]{\textbf{Technique 6:}\textit{
		Every arc connecting a node $i$ to itself, e.g., $(i,i)$, can be eliminated.}
	
	Regardless of whether node $i$ is an origin, transshipment, or destination node, it is sub-optimal for the residents associated with an origin to visit a node twice. By this, all of the loops can be eliminated from $\mathbb{G}$ without affecting the remaining network (Figure \ref{tech}i). }

\end{itemize}
Technique 7 is a pruning technique called node elimination and valency proposed by \cite{ogbuobiri1970sparsity}.
\begin{itemize}
\item[]{\textbf{Technique 7:}\textit{
		Assume there exists a non-vulnerable path $i\rightarrow j \rightarrow k$ (or with the opposite direction), where $j \in \mathcal{N}_\tau$, and nodes $i,k \in \mathcal{N}$ are the only neighbors of node $j$. Then, eliminating node $j$ and the arcs incident to it and adding a new arc between $i$ and $k$ does not interrupt the flow in $\mathbb{G}$ (Figure \ref{tech}j). }
	
	Since node $j$ is a transshipment node, substituting the path  $i\rightarrow j \rightarrow k$ with a direct arc from $i$ to $j$, i.e., $(i,j)$ will not change the flow in $\mathbb{G}$. Consequently, $t_{ik}=t_{ij}+t_{jk}$ (or $t_{ki}=t_{kj}+t_{ji}$) represents the travel time of the newly added arc. In a special case where $i=k$, node $j$ and arcs $(i,j)$ and $(j,k)$ can be eliminated.
}

\end{itemize}
Finally, in Technique 8, we introduce a new arc pruning approach that exploits the properties of graph cliques. The implementation to discover the cliques has been adopted from an algorithm developed by \cite{zhang2005genome}, taking advantage of large-memory architectures and fixed-parameter tractability that reduce the search space and produce computationally feasible performance.
\begin{itemize}
\item[]{\textbf{Technique 8:}\textit{
		By finding the existing cliques of size three in $\mathbb{G}$ containing non-vulnerable arcs $(i,j)$, $(i,h)$, and $(j,h)$, it is possible to remove arc $(i,h)$ if $t_{ih}\ge t_{ij}+t_{jh}$ (Figure \ref{tech}k). }
	
	The explanation is similar to the properties of the triangle inequality. The difference is that there the sign of inequality is reversed. Otherwise, if $t_{ih}< t_{ij}+t_{jh}$, the following valid inequality can be added to RNFMP:
	\begin{align}		 
		&x_{ij}^k+x_{ih}^k+x_{jh}^k \le 1 &\forall k \in \mathcal{N}_o, \forall i,j,h \in \mathcal{C}\label{VI2}
	\end{align}
	where $\mathcal{C}$ denotes the set of non-vulnerable cliques of size three in $\mathbb{G}$.		
}
\end{itemize}
\vspace{-10pt}
\subsection{Valid inequalities and variable elimination techniques}\label{VIVR}	
The techniques introduced in Section \ref{PT} consist of network pruning methods that will reduce the number of decision variables and size of the search space. In this section, the proposed techniques will reduce the number of decision variables without further pruning the network. 

In Proposition \ref{yVI}, a valid inequality is introduced along with a decision variable elimination technique focusing on eliminating redundant $y_{ij}$ variables. 
	\begin{proposition}\label{yVI}
		If there exists a node $i\in \mathcal{N}_o$ that is not incident to any non-vulnerable outgoing arc (i.e., $\nexists j \in \mathcal{N}: (i,j) \in \mathcal{A}_n$), and $\delta_{out}(i)\ge 1$, then:
		\begin{align}		 
			&\sum_{j:(i,j) \in \mathcal{A}_v} y_{ij} \ge 1 \label{VI}
		\end{align}
	\end{proposition}
The proof is given in the e-companion to this paper (Section \ref{proofs}).

In Proposition \ref{prop2}, we introduce a new idea, developed based on the unique characteristics of RNFMP, to decrease the number of $x_{ij}^k$ variables.
\begin{proposition}\label{prop2}
For every origin $k\in \mathcal{N}_o$, which is connected to all of the existing destinations in the non-upgraded network  $(i.e., \mathbb{G}\setminus \mathcal{A}_v)$, variable $x_{ij}^k$ can be eliminated from RNFMP if
\begin{align}		 
	&SP_f(k,i)+t_{ij}> \max_{g \in \mathcal{N}_d}\{SP_n(k,g)\} & \forall (i,j)\in \mathcal{A} \label{p2}
\end{align}
where $SP_n(u,v)$ and $SP_f(u,v)$ are the lengths of the shortest paths from node $u\in \mathcal{N}_o$ to node $v\in \mathcal{N}$ in the non-upgraded $(\mathbb{G}\setminus \mathcal{A}_v)$ and fully upgraded $(\mathbb{G})$, respectively.
\end{proposition}
The proof is given in the e-companion to this paper (Section \ref{proofs}).

Unlike Proposition \ref{yVI}, to implement the variable elimination technique explained in Proposition \ref{prop2}, we need to perform a polynomial-time algorithm such as Dijkstra's algorithm, $|\mathcal{N}_d|$ times for each origin $k$, to find the shortest path from $k$ to every $|\mathcal{N}_d|$ destinations in $\mathbb{G}\setminus \mathcal{A}_v$. The length of the longest shortest path among them will be selected to be compared with the length of the shortest paths from $k$ to every arc $(i,j)$ in the upgraded network $\mathbb{G}$. The $x_{ij}^k$ variables that correspond to the $(i,j)$ arcs satisfying Inequality \eqref{p2} have been kept in the model, and the rest will be discarded.

In the next proposition, similar to pruning Technique 1, we will leverage the APs.
\begin{proposition}\label{delta}
If there exists a disconnected component $\Delta$ obtained from $\mathbb{G}\setminus\{v\}$, where $v\in \mathcal{N}_{p}$ and $\nexists u\in \mathcal{N}_{\Delta} : u \in \mathcal{N}_d$, then $x_{ij}^k$ variables can be eliminated for every
$k \notin \mathcal{N}_{\Delta}$ and $(i,j)\in\mathcal{A}_{\Delta}$.
\end{proposition}
The proof is given in the e-companion to this paper (Section \ref{proofs}).
\subsection{Initial solution}\label{init}
Creating a good initial solution can often help an integer programming solver find the optimal solution faster. Consequently, here we propose a greedy approach to determine the initial solution for RNFMP. In this approach, each origin $k$ is associated with two vectors, $D$ and $D^\prime$. Both vectors represent a permutation of destinations sequenced in order of the closest to furthest in terms of travel time from $k$. For both of these vectors, all travel time values are calculated based on the shortest paths between the OD pairs. For $D^\prime$, the travel times between $k$ and destinations have been calculated excluding vulnerable roads in $\mathbb{G}$. Accordingly, the size of vector $D$ is always equal to the number of destinations (i.e., $|\mathcal{N}_d|$). However, since removing vulnerable roads may disrupt the connectivity of the network, we have $0\le |D^\prime| \le |\mathcal{N}_d|$. For instance, if $|\mathcal{N}_d|=3$, then $D^\prime=\{D_3,D_1\}$ indicates that the closest and second closest destinations to origin $k$ are $D_3$ and $D_1$, respectively. However, $k$ has no connected path to the second destination ($D_2$). Table \ref{IFS} contains some examples for these vectors in the fourth column ($D$ \& $D^\prime$) for different origins.

To generate the initial solution, we sort the set of origins based on their population size. Starting with the most populated origin $k=\argmaxB_{i\in \mathcal{N}_o} \{h^i\}$, at each iteration, the origin will be assigned to its closest destination in $D^\prime$. However, before conducting the assignment, the feasibility of Constraints \eqref{8} are checked by comparing the population of the origin $k$ and the remaining capacity at the closest destination. This comparison is shown in the fifth column of Table \ref{IFS} (HCFs' Availability).

In the example provided in Table \ref{IFS}, $k$ is the most populated origin followed by origins $l$, $g$, and $m$, and $D^\prime=\{D_3,D_1,D_2\}$. This implies that, at iteration 1, origin $k$ should be assigned to the third destination ($O_{k} \rightarrow D_3$). At the end of each iteration, we update the destination's available capacity based on the population of origin $k$. In the provided example, at iterations 1 and 2, all of the destinations are available. Thus, the closest destination, which is $D_3$, is selected, and its capacity is reduced by $h^k$ and $h^l$ in iterations 1 and 2, respectively. At each iteration, if $D^\prime = \O$, then the associated origin will be assigned to its closest destination based on vector $D$. In the third iteration, the third-largest origin ($O_{g}$) does not have access to any destination through non-vulnerable roads. Therefore, we determine its closest destination based on vector $D$, which is $D_3$. However, due to assigning the previous origin to $D_3$, its capacity is no longer enough to serve $O_{g}$ with $h^g$ residents. As a result, this origin will be assigned to $D_1$, the next available destination in $D$.

\vspace{-10pt}
\begin{table}[!ht]
\caption{An example for the initial solution approach}
	\centering
		\tiny
	{\def\arraystretch{1}  
		\begin{tabular}{cccllll}
		\hline\hline
		Iteration & Origin & Population & $D$ \& $D^\prime$ & HCFs' Availability & OD Assignment & Capacity Update \\
		\hline
\multirow{1}[2]{*}{1} & \multirow{1}[2]{*}{k} & \multirow{1}[2]{*}{$h^k$} & $D^\prime=\{D_3,D_1,D_2\}$ & $H_1>h^k:$ \checkmark & \multirow{1}[2]{*}{$O_{k} \rightarrow D_3$} & \multirow{1}[2]{*}{$H_3=H_3-h^k$} \\
&       &       & \multirow{2}[1]{*}{$D=\{D_3,D_1,D_2\}$} & $H_2>h^k:$ \checkmark &       &  \\
&       &       &       & $H_3>h^k:$ \checkmark &       &  \\
		\hline
\multirow{1}[2]{*}{2} & \multirow{1}[2]{*}{$l$} & \multirow{1}[2]{*}{$h^l$} & $D^\prime=\{D_3,D_2\}$ & $H_1>h^l:$ \checkmark & \multirow{1}[2]{*}{$O_{l} \rightarrow D_3$} & \multirow{1}[2]{*}{$H_3=H_3-h^l$} \\
&       &       & \multirow{2}[1]{*}{$D=\{D_1,D_3,D_2\}$} & $H_2>h^l:$ \checkmark &       &  \\
&       &       &       & $H_3>h^l:$ \checkmark &       &  \\
		\hline
\multirow{1}[2]{*}{3} & \multirow{1}[2]{*}{$g$} & \multirow{1}[2]{*}{$h^g$} & $D^\prime=\{\}$ & $H_1>h^g:$ \checkmark & \multirow{1}[2]{*}{$O_{g} \rightarrow D_1$} & \multirow{1}[2]{*}{$H_1=H_1-h^g$} \\
&       &       & \multirow{2}[1]{*}{$D=\{D_3,D_1,D_2\}$} & $H_2>h^g:$ \checkmark &       &  \\
&       &       &       & $H_3>h^g :$ X &       &  \\
		\hline
\multirow{1}[1]{*}{4} & \multirow{1}[1]{*}{m} & \multirow{1}[2]{*}{$h^m$}   & $D^\prime=\{D_1\}$ & $H_1>h^m:$ \checkmark & \multirow{1}[1]{*}{$O_{m} \rightarrow D_1$} & \multirow{1}[1]{*}{$H_1=H_1-h^m$} \\
&       &       & \multirow{2}[0]{*}{$D=\{D_1,D_2,D_3\}$} & $H_2>h^m:$ \checkmark &       &  \\
&       &       &       & $H_3>h^m:$ X &       &  \\
\multirow{2}[0]{*}{\begin{turn}{90}...\end{turn}} & \multirow{2}[0]{*}{} & \multirow{2}[0]{*}{} & \multirow{2}[0]{*}{} & \multirow{2}[0]{*}{} & \multirow{2}[0]{*}{} & \multirow{2}[0]{*}{} \\
&       &       &       &       &       &  \\
$|\mathcal{N}_d|$ &       &       &       &       &       &  \\
			\hline\hline
		\end{tabular}
	\label{IFS}%
	}
\end{table}
\vspace{-20pt}

	\section{Data collection and design of experiments}\label{DOE}
For our experiments, we consider the cities of Fort Dodge and Coralville in the state of Iowa. Iowa has experienced two major flood events in recent years (1993 and 2008), so road mitigation is an important topic in the state \citep{zogg2014top,DOT2019}. We use these two cities as samples of road networks of different sizes with different numbers of HCFs.  We will computationally explore the performance of our solution techniques on both cities in Section \ref{impimp}.  In Section \ref{results}, we will dive deeper into the implications of the solutions for Coralville, using it as an in-depth case study.   

\subsection{Instance size}
For both cities, the road network is available from the \cite{DOT2019} and includes the roads within city limits plus a three-mile radius. Using the ArcGIS program, we identify the length of each road. Each road's travel time is then calculated by dividing the length by its speed limit. The network of Fort Dodge, which represents a smaller road network than Coralville, consists of 1,688 nodes formed as road intersections and 4,752 directed arcs (road links). For Coralville, it consists of 3,841 nodes and 9,611 directed arcs.

\subsection{Vulnerable roads}  
\begin{figure}[h!]
\begin{center}
	\caption{Coralville roads within a 100-year floodplain\\}
	\vspace{-10pt}
\includegraphics[scale=0.4]{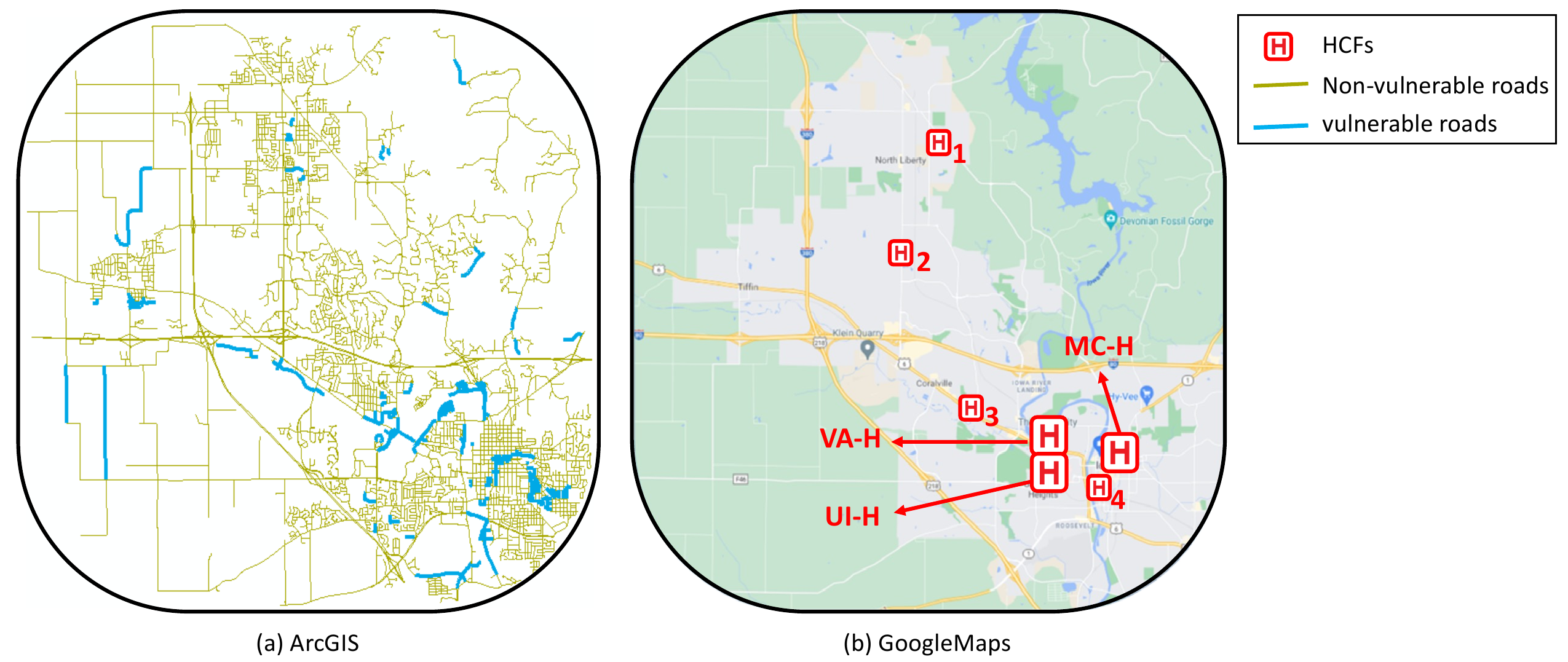}
\label{flooded}
\end{center}
\end{figure}

Two representations of the road network obtained from the Iowa Department of Transportation, generated by two GIS programs, ArcGIS and GoogleMaps, are shown in Figure \ref{flooded} for the city of Coralville. 
In Figure \ref{flooded}a, thick blue lines show the road segments that lie inside the 100-year floodplain, the \textit{vulnerable} roads in RNFMP.  By comparing this map with the one generated with GoogleMaps in Figure \ref{flooded}b, we can see that the majority of vulnerable roads are located at the southeast of the network where the Iowa river passes through the residential area. The occurrence of this type of flood in the region would result in the closure (flooding) of 268 road segments (472 directed arcs). The type of all existing roads and the number of vulnerable roads of each type are illustrated in Table \ref{types} \citep{wiki:xxx,ICE}. 
Following the same approach, 140 vulnerable arcs were identified for Fort Dodge. 
\vspace{-10pt}
\vspace{-35pt}
{\footnotesize
\begin{table}[h!]
\centering
\tiny
\caption{The types of roads in our network}
\renewcommand{\arraystretch}{1}
\begin{tabular}{lrrp{30.145em}}
	\hline\hline
	Type  & {\# roads} & {\# vulnerable} & \multicolumn{1}{l}{Description} \\
	\hline
	Motorway &                  123  & 0     & Major divided highways with emergency shoulder (e.g., freeway) \\
	Primary &                  198  & 14    & \multicolumn{1}{l}{State or major highways (link large towns)} \\
	Secondary &                  461  & 32    & \multicolumn{1}{l}{County highways (link towns)} \\
	Tertiary &                  496  & 25    & \multicolumn{1}{l}{Roads with low to moderate traffic (link small towns and villages)} \\
	Unclassified &                  133  & 5     & \multicolumn{1}{l}{Local or minor roads (often link villages)} \\
	Residential &               3,913  & 188   & Roads with the lowest speed limits and capacities (provide access to housing or within residential areas) \\
	Rest area &                       4  & 0     & Place where drivers can leave the road to rest \\
	\hline\hline
\end{tabular}%
\label{types}%
\end{table}%
}
\vspace{-25pt}

\subsection{Population centers}
To obtain the location of population centers, we have used the demographic data, including a set of centroids, received from the Iowa DOT. Each centroid represents the population for a given area distinguished by a city block. Using ArcGIS, the population of each block has been divided by the number of nodes within the block, and equally distributed among them. This gives us an estimate of the number of residents associated with each individual node. In the experiments, the number of residents associated with each origin is used as its weight in RNFMP's objective function.
\vspace{-10pt}
\begin{table}[htbp]
	\tiny
	\centering
	\caption{Demographic features of each class of instances}
	\renewcommand{\arraystretch}{1}
	\begin{tabular}{llcccc}
		\hline\hline
		&       & \textit{p} & \# residents & \# population centers ($|\mathcal{N}_o|$) & Covered population \\
		\hline
		Class 1:     &       & 24    & 62,721 & \mbox{ }\mbox{ }859   & 80\% \\
		Class 2:     &       & 18    & 66,863 & 1,062 & 85\% \\
		Class 3:     &       & 14    & 70,560 & 1,296 & 90\% \\
		Class 4:     &       & 9     & 74,521 & 1,666 & 95\% \\
		\hline\hline
	\end{tabular}%
	\label{pop3}%
\end{table}%
\vspace{-10pt}
In each generated instance, the set of population centers ($\mathcal{N}_o$) will be generated based on the minimum number of residents $p$ associated with each node $i\in \mathcal{N}$. A node qualifies as a population center if the number of residents associated with the node is not less than $p$. Decreasing this value would increase the cardinality of the set $\mathcal{N}_o$ and thus would increase the computational cost. 
For instance, in Table \ref{pop3}, by setting $p=24$ or $p=9$ the number of residents (i.e., $\sum_{k \in \mathcal{N}_o}h^k$) that will be covered equals 62,721 and 74,521, which are about 80\% and 95\% of Coralville's total population, respectively.
Four classes of instances have been carried out in our experiments, considering 80\%, 85\%, 90\%, and 95\% of the total population. Table \ref{pop3} presents the $p$ level values and the corresponding populations for Coralville for each of the four classes. For Fort Dodge, the corresponding $p$ values for classes 1 through 4 are equal to 16, 13, 10, and 6, respectively.

\subsection{HCFs and capacity levels}

\begin{figure}[h]
\begin{center}
		\caption{Testbed networks}
	\includegraphics[scale=0.3]{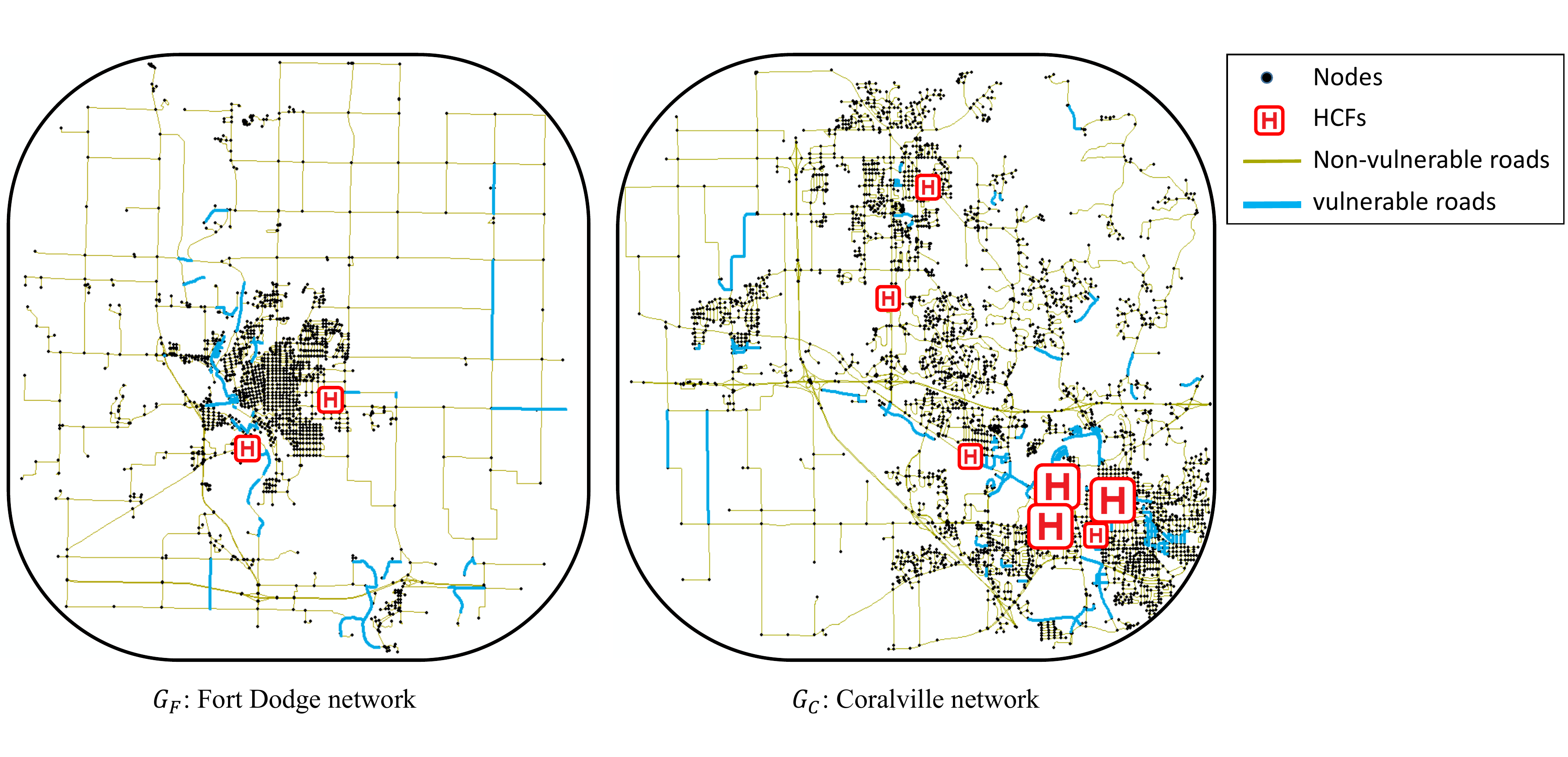}
	\label{testbed}
\end{center}
\end{figure}

Based on the \cite{DOT2019},  there are two healthcare facilities in Fort Dodge. These facilities are set as the networks' destinations and are labeled with the \vahid{\mybox{\textbf{H}}} signs in Figure \ref{testbed}. There exist seven HCFs in the network for Coralville as shown in Figure \ref{flooded}b and Figure \ref{testbed}. Among these facilities, only three of them are hospitals (shown by larger \vahid{\mybox{\textbf{H}}} signs), and the rest are smaller facilities providing mostly non-emergency or routine care (e.g., clinics).  We first set all of these seven facilities as the set of destinations. Next, we consider the set of destinations based on only the set of hospitals: i) \textit{University of Iowa Hospitals \& Clinics (UI-H)}, ii) \textit{Mercy Iowa City (MC-H)},  iii) \textit{Iowa City VA Health Care System (VA-H)}, as shown in Figure \ref{flooded}a.
\vspace{-20pt}
\vspace{-20pt}

An identical capacity  is used for each healthcare facility (i.e., $H_j=H$, $\forall j \in \mathcal{N}_d$) in which the total capacity of all HCFs is considered $\alpha$ percent higher than the total population in Fort Dodge and one set of experiments for Coralville. The capacity of each hospital is defined as $H={\dfrac{(1+\alpha)\sum_{k \in \mathcal{N}_o}h_k}{|\mathcal{N}_d|}}$, where $\alpha=\{15\%,30\%\}$.

We also consider a set of instances for the city of Coralville where the capacity of each hospital is calculated based on the proportion of the actual number of beds available at each hospital. As of January 2022, there were 796, 183, and 93 beds available in UI-H, MC-H, and VA-H, respectively \citep{UI-H,MC-H,VA-H}. This indicates that, in case of considering these three hospitals as the only available HCFs in the network, around 75\% of the total capacity will be allocated to UI-H.

\subsection{Upgrading cost}\label{cost}
There are two common approaches to measure road upgrading costs in the literature.  The cost to upgrade each road is either (i) equal to one unit, where the total budget is defined as the total number of upgradeable roads (e.g., see \cite{angulo2016lagrangian} and \cite{bagloee2017identifying}), or (ii) it is calculated as a weighted \textit{score} proportional to (as a linear function of) its length (e.g., see \cite{peeta2010pre} and \cite{duque2011grasp}). The current paper, unlike the literature, considers actual upgrading costs (\$) for roads and bridges. 

Based on the road construction costs reported by the Midwest Economic Policy Institute, the upgrading cost in the State of Iowa per mile per lane is roughly equal to \$32,000 \citep{craighead2018comparison}. This includes construction (\$26,436), right-of-way acquisition (\$1,684), and  engineering (\$3,977) costs. 

\subsection{Mitigation budget}\label{mitbud}
We define different scenarios based on various budget values set as the percentage of the total investment required to upgrade all vulnerable road segments ($\hat{B}$). In Section \ref{sensb}, we will show that a mitigation budget of \$6,001,287, which is about 15\% of $\hat{B}$ in our original network, is enough to achieve the same efficiency level of our fully functional network. Additionally, any amount of budget less than 5\% of $\hat{B}$ would result in a disconnected network (infeasibility) in a subset of instances for Coralville. As a result, four different scenarios based on four budget levels of 5\%, 7.5\%, 10\%, and 15\% are defined for instances with identical capacities. Following the same approach, three budget levels of 7.5\%, 10\%, and 15\% are defined for Fort Dodge. For our experiments in Section \ref{results}, we generate $4\times4\times4\times2=128$ instances of the RNFMP for Coralville based on combinations of (budget, HCFs' total capacity, covered percentage of the total population, and the number of HCFs), where budget $\in \{5\%, 7.5\%, 10\%, 15\%\}$ of $\hat{B}$, HCFs' total capacity is $\alpha \in \{15\%, 30\%\}$ higher than the total population (for instances with either identical and different capacities based on the number of beds), covered percentage of the total population $\in\{80\%, 85\%, 90\%, 95\%\}$, and number of HCFs $\in \{3,7\}$.

\section{Evaluation of model improvements}\label{impimp}
In this section, we evaluate the effectiveness of the ideas presented in Section \ref{imp3} for Fort Dodge and Coralville.
For the test instances, we use a subset of 16 instances out of these 128 instances for Coralville plus a set of 12 instances of the RNFMP for Fort Dodge. The features of these 28 instances are summarized in Table \ref{inst}. The second column specifies the $p$ values to cover 80, 85, 90, and 95 percent of the population. In the last column,
it can be seen that in total, 12 and 16 instances are created for Fort Dodge ($\mathbb{G}_F$) and Coralville ($\mathbb{G}_C$) networks, respectively.

\begin{table}[htbp]
	\tiny
	\centering
	\caption{Characteristics of test instance}
		\renewcommand{\arraystretch}{1}
	\begin{adjustbox}{max width=\textwidth}	
		\begin{tabular}{lllll}
			\hline\hline
			Network & $p$ & \# HCFs & Budget & \# instances \\
			\hline
			$\mathbb{G}_F$ & \{16,13,10,6\} & \{2\} & \{7.5\%,10\%,15\%\} & 12 \\				
			$\mathbb{G}_C$ & \{24,18,14,9\} & \{3\} & \{5\%,7.5\%,10\%,15\%\} & 16 \\
			\hline\hline
		\end{tabular}%
	\end{adjustbox}
	\label{inst}%
	\label{testinst}
\end{table}%

In Section \ref{effprun}, we evaluate the impact of pruning techniques in terms of the number of eliminated nodes, arcs, and variables. The variable elimination methods are evaluated with regard to the percentage of variables that are eliminated. Finally, in Section \ref{imp_results}, 
we verify the impact of all of the improvements, including the creation of an initial feasible solution and valid inequalities, on the CPU runtime.  All of our pruning and variable elimination techniques are implemented in Python 3.7.0 programming language. We use Gurobi version 9.1 provided by Gurobi Optimization, Inc. (2021) to solve the experiments. The IP is solved with a 32 thread count on the University of Iowa's Argon high-performance computing cluster \citep{argon}.

\subsection{Problem size reduction}\label{effprun}
To verify the impact of the pruning techniques, we compare the number of nodes, arcs, and decision variables in each network before and after applying the different techniques. Note that changing the budget does not change the \textit{size} of the network, so we do not consider different budgets for instances created in Table \ref{testinst}.

Table \ref{elims} shows the results for each network and instance, specified in the first two columns. The next two columns show the parameter setting for each instance, and the fifth column shows the number of variables before applying the pruning techniques. Finally, the last three columns show the percentage of variables, nodes, and arcs eliminated, respectively.
The average percent of eliminated variables for $\mathbb{G}_F$ and $\mathbb{G}_C$ is equal to 53.17\% and 54.29\%, respectively. 
\vspace{-10pt}
\begin{table}[htbp]
	\tiny
	\renewcommand{\arraystretch}{1}
	\centering
	\caption{Percentage of variable, arc, and node eliminations for each instance}	
	\begin{adjustbox}{max width=\textwidth}	
		\begin{tabular}{llllrrrr}
			\hline\hline
			Network & \multicolumn{1}{c}{Instance} & \multicolumn{1}{c}{$p$} & \multicolumn{1}{c}{HCFs} & \multicolumn{1}{c}{Original } & \multicolumn{3}{c}{Eliminations} \\
			\cmidrule{6-8}          &       &       &       & \# Variables & Variables & Nodes & Arcs \\
			\hline
			$\mathbb{G}_F$ & \multicolumn{1}{c}{1} & \multicolumn{1}{r}{16} & \multicolumn{1}{c}{2} & 2,399,900 & 54.30\% & 20.79\% & 18.83\% \\
			& \multicolumn{1}{c}{2} & \multicolumn{1}{r}{13} & \multicolumn{1}{c}{2} & 2,917,868 & 53.20\% & 19.31\% & 17.70\% \\
			& \multicolumn{1}{c}{3} & \multicolumn{1}{r}{10} & \multicolumn{1}{c}{2} & 3,416,828 & 53.08\% & 18.66\% & 17.23\% \\
			& \multicolumn{1}{c}{4} & \multicolumn{1}{r}{6} & \multicolumn{1}{c}{2} & 4,295,948 & 52.09\% & 15.40\% & 14.88\% \\
			\cmidrule{5-8}          &       &       &       & Avg. (\%) = & 53.17\% & 18.54\% & 17.16\% \\
			\hline
			$\mathbb{G}_C$ & \multicolumn{1}{c}{5} & \multicolumn{1}{r}{24} & \multicolumn{1}{c}{3} & 8,256,321 & 57.06\% & 26.37\% & 23.27\% \\
			& \multicolumn{1}{c}{6} & \multicolumn{1}{r}{18} & \multicolumn{1}{c}{3} & 10,165,739 & 55.72\% & 24.11\% & 21.42\% \\
			& \multicolumn{1}{c}{7} & \multicolumn{1}{r}{14} & \multicolumn{1}{c}{3} & 12,376,149 & 54.73\% & 21.95\% & 19.54\% \\
			& \multicolumn{1}{c}{8} & \multicolumn{1}{r}{9} & \multicolumn{1}{c}{3} & 15,884,587 & 54.15\% & 19.21\% & 17.11\% \\
			\cmidrule{5-8}          &       &       &       & Avg. (\%) = & 55.41\% & 22.91\% & 20.34\% \\
			\hline
			Avg. (\%) = &       &       &       &       & 54.29\% & 20.73\% & 18.75\% \\
			\hline\hline
		\end{tabular}%
	\end{adjustbox}
	\label{elims}
\end{table}

\vspace{-15pt}
Table \ref{sizered2} demonstrates the problem size reduction by technique for Instance 5 from Table \ref{elims}. The table's first row shows the network's size ($\mathbb{G}_C$) in terms of the number of variables, nodes, and arcs. The next 11 rows show the impact of each method individually in reducing the size of the problem. The first column specifies the applied technique, the second column reports the reduction in the number of RNFMP's \textit{decision variables}, and in the last two columns, we report the reductions in the network's number of \textit{nodes} and \textit{arcs}, respectively. In each column, the percentage of eliminated items is represented in the `` (\%)" sub-columns. Pruning Techniques 1 and 7, along with Proposition \ref{delta}, are the most powerful methods in terms of their impact on the number of eliminated decision variables. The impact of the techniques in other instances is almost the same, indicating that not only Techniques 1 and 7, and Proposition \ref{delta} are the most effective methods in this instance, but the most effective ones in general. The last row in Table \ref{sizered2} shows the reductions after applying all of the improvements together. For the specific instance highlighted in Table \ref{sizered2}, the proposed techniques reduce the number of decision variables, nodes, and arcs by almost 34\%, 26\%, and 25\%, respectively.
\vspace{-10pt}
\begin{table}[!htb]
	\tiny
	\renewcommand{\arraystretch}{1}
	\centering
	\caption{Percentage of variable, arc, and node eliminations for Instance 9}	
	\begin{adjustbox}{max width=\textwidth}	
			\begin{tabular}{lrrrrrrrrr}
				\hline\hline
				&       & \multicolumn{2}{c}{Variables} &       & \multicolumn{2}{c}{Nodes} &       & \multicolumn{2}{c}{Arcs} \\
				\cmidrule{3-4}\cmidrule{6-7}\cmidrule{9-10}     Size of the original network: &       & 8,256,321 & (\%)  &       & 3,841 & (\%)  &       & 9,611 & (\%) \\
				\cmidrule{3-4}\cmidrule{6-7}\cmidrule{9-10}    Technique 1 &       & -1,577,604 & 19.11\% &       & -917  & 23.87\% &       & -1,837 & 19.11\% \\
				Technique 2 &       & -20,144 & 0.24\% &       & -16   & 0.42\% &       & -24   & 0.25\% \\
				Technique 3 &       & -114,634 & 1.39\% &       & -67   & 1.74\% &       & -153  & 1.59\% \\
				Technique 4 &       & -16,708 & 0.20\% &       & -6    & 0.16\% &       & -24   & 0.25\% \\
				Technique 5 &       & -176,567 & 2.14\% &       & N/A   & N/A   &       & -205  & 2.13\% \\
				Technique 6 &       & -57,944 & 0.70\% &       & N/A   & N/A   &       & -68   & 0.71\% \\
				Technique 7 &       & -1,218,449 & 14.76\% &       & -928  & 24.16\% &       & -1,850 & 19.25\% \\
				Technique 8 &       & -43,339 & 0.52\% &       & N/A   & N/A   &       & -51   & 0.53\% \\
				Proposition 1 &       & -462  & 0.01\% &       & N/A   & N/A   &       & N/A   & N/A \\
				Proposition 2 &       & -481,401 & 5.83\% &       & N/A   & N/A   &       & N/A   & N/A \\
				Proposition 3 &       & -889,850 & 10.78\% &       & N/A   & N/A   &       & N/A   & N/A \\
				\hline
				All improvements &       & -2,789,067 & 33.78\% &       & -1,012 & 26.3\% &       & -2,392 & 24.89\% \\
				\hline\hline
			\end{tabular}%
		\end{adjustbox}
		\label{sizered2}
	\end{table}

	\vspace{-10pt}
	Note that covering a larger percentage of the population requires a larger set of origins ($|\mathcal{N}_o|$). With more origins, more nodes and arcs need to be traversed on the network, which prevents pruning and eliminating as many nodes and arcs. Similarly, increasing the number of destinations lowers the percentage of eliminations.

	\subsection{Impact on solution time}\label{imp_results}
	To verify the impact of the improvements on the CPU runtime, we first solve the RNFMP model with the original network ($\mathbb{G}$). Next, we compare the runtime with those of the pruned network, including the use of the initial solution and valid inequalities. The runtime to implement the developed improvements is negligible. Preliminary results show that even in large instances with a solution time longer than three hours, implementing all of the improvements takes less than six minutes. Thus we only focus on reporting the model's solution time.
	
	Table \ref{imp-FD} demonstrates the quality of improvements for the Fort Dodge network.
	The first three columns show the parameter setting for each instance. In the remaining columns, we report the solution time of RNFMP for each instance before and after adding different improvements. We refer to the pruned network as $\mathbb{G}^\prime$. Accordingly, $\mathbb{G}^\prime_{IS}$ and $\mathbb{G}^\prime_{VI}$ columns report the runtime of the RNFMP model with the pruned network after adding the improvements with an initial solution (IS) and valid inequalities (VI), respectively. The runtime after adding all of the improvements is reported in the last column, $\mathbb{G}^\prime_{IS+VI}$. The instances are solved with a three-hour time-limit (TL), and in each column, runtime (RT) values reported are in seconds. The percentage of runtime improvement over the value in column 3 is presented in the ``Imp" sub-columns.
	\begin{table}[!htb]
		\tiny
		\renewcommand{\arraystretch}{1}
		\centering
		\caption{Impact of the improvements for Fort Dodge}
		\begin{adjustbox}{max width=\textwidth}	
			\begin{tabular}{rlrrrrrrrrrrrrrr}
				\hline\hline
				\multicolumn{1}{l}{Instance} & \multicolumn{1}{c}{Budget} & \multicolumn{1}{c}{$p$} & \multicolumn{1}{c}{$\mathbb{G}$} &       & \multicolumn{2}{c}{$\mathbb{G}^\prime$} &       & \multicolumn{2}{c}{$\mathbb{G}^\prime_{IS}$} &       & \multicolumn{2}{c}{$\mathbb{G}^\prime_{VI}$} &       & \multicolumn{2}{c}{$\mathbb{G}^\prime_{IS+VI}$} \\
				\cmidrule{6-7}\cmidrule{9-10}\cmidrule{12-13}\cmidrule{15-16}          &       &       & \multicolumn{1}{c}{RT} &       & \multicolumn{1}{c}{RT} & \multicolumn{1}{c}{Imp} &       & \multicolumn{1}{c}{RT} & \multicolumn{1}{c}{Imp} &       & \multicolumn{1}{c}{RT} & \multicolumn{1}{c}{Imp} &       & \multicolumn{1}{c}{RT} & \multicolumn{1}{c}{Imp} \\
				\hline
				\multicolumn{1}{l}{1} & \multicolumn{1}{c}{15\%} & \multicolumn{1}{c}{16} & 370.73 &       & 152.75 & 59\%  &       & 159.87 & 57\%  &       & 136.12 & 63\%  &       & 145.84 & 61\% \\
				\multicolumn{1}{l}{2} &       & \multicolumn{1}{c}{13} & 420.50 &       & 265.27 & 37\%  &       & 192.75 & 54\%  &       & 179.61 & 57\%  &       & 235.75 & 44\% \\
				\multicolumn{1}{l}{3} &       & \multicolumn{1}{c}{10} & 845.29 &       & 450.94 & 47\%  &       & 416.98 & 51\%  &       & 285.48 & 66\%  &       & 254.98 & 70\% \\
				\multicolumn{1}{l}{4} &       & \multicolumn{1}{c}{6} & 842.95 &       & 450.85 & 47\%  &       & 525.82 & 38\%  &       & 530.26 & 37\%  &       & 315.72 & 63\% \\
				\multicolumn{1}{l}{5} & \multicolumn{1}{c}{10\%} & \multicolumn{1}{c}{16} & 371.46 &       & 165.81 & 55\%  &       & 164.31 & 56\%  &       & 182.06 & 51\%  &       & 161.96 & 56\% \\
				\multicolumn{1}{l}{6} &       & \multicolumn{1}{c}{13} & 1,331.79 &       & 686.86 & 48\%  &       & 712.65 & 46\%  &       & 770.60 & 42\%  &       & 673.80 & 49\% \\
				\multicolumn{1}{l}{7} &       & \multicolumn{1}{c}{10} & 1,594.85 &       & 827.70 & 48\%  &       & 758.94 & 52\%  &       & 938.16 & 41\%  &       & 754.11 & 53\% \\
				\multicolumn{1}{l}{8} &       & \multicolumn{1}{c}{6} & 965.66 &       & 371.28 & 62\%  &       & 389.99 & 60\%  &       & 419.19 & 57\%  &       & 371.46 & 62\% \\
				\hline
				\multicolumn{1}{l}{9} & \multicolumn{1}{c}{7.50\%} & \multicolumn{1}{c}{16} & 386.56 &       & 162.35 & 58\%  &       & 164.43 & 57\%  &       & 146.80 & 62\%  &       & 168.19 & 56\% \\
				\multicolumn{1}{l}{10} &       & \multicolumn{1}{c}{13} & 4,686.63 &       & 743.47 & 84\%  &       & 722.28 & 85\%  &       & 745.08 & 84\%  &       & 683.01 & 85\% \\
				\multicolumn{1}{l}{11} &       & \multicolumn{1}{c}{10} & 1,740.82 &       & 928.60 & 47\%  &       & 723.35 & 58\%  &       & 698.27 & 60\%  &       & 896.82 & 48\% \\
				\multicolumn{1}{l}{12} &       & \multicolumn{1}{c}{6} & 791.89 &       & 344.47 & 57\%  &       & 355.80 & 55\%  &       & 333.31 & 58\%  &       & 383.83 & 52\% \\
				\hline
				& Average: &       & 1,195.76 &       & 462.53 & 54\%  &       & 440.60 & 56\%  &       & 447.08 & 57\%  &       & 420.46 & 58\% \\
			\end{tabular}%
		\end{adjustbox}
		\label{imp-FD}%
	\end{table}%

	After pruning the network, the average runtime for  $\mathbb{G}^\prime$ improves by 54\% compared to the original Fort Dodge network. The average runtime for  $\mathbb{G}^\prime_{IS}$ and $\mathbb{G}^\prime_{VI}$ is improved by 56\% and 57\% across the 12 instances. After applying all of the improvements, the average runtime of the instance reduces by 58\% in $\mathbb{G}^\prime_{IS+VI}$. This indicates that even though providing an initial solution and adding the valid inequalities elevate the average quality of the improvements, the main contribution to the runtime reductions for Fort Dodge belongs to the developed pruning techniques (Techniques 1 through 8).
	
	In Section \ref{impact-coral}, Table \ref{imp-C}, similar to $\mathbb{G}_F$, the quality of the improvements across the 16 instances generated for the $\mathbb{G}_C$ network is reported. For the instances that Gurobi was not able to solve to optimality within the time-limit, we report the MIPGap at termination inside the parenthesis, i.e., TL (MIPGap\%). This table highlights the significant impact of the improvements. For example, 
	while the MIPGap in instances 10 and 12 equal to 0.01\% and 0.03\% (column $\mathbb{G}$), after applying all of the improvements (column $\mathbb{G}^\prime_{IS+VI}$), Gurobi could solve them in 4,521 and 4,791 seconds, respectively.

	\section{Case Study}\label{results}
In the following, we examine the results for the city of Coralville in more detail. In this section, we look at the effect each parameter has on the selected upgraded roads and the objective to help provide managerial insights for the city of Coralville for mitigation strategies.

\subsection{Sensitivity analysis of mitigation budget}\label{sensb}
For each combination of $p$, $\alpha$, $|\mathcal{N}_d|$, and the value of the HCF's capacity, illustrated in Section \ref{mitbud}, we  solve the RNFMP over the fully functional network where the budget constraint is relaxed and all vulnerable arcs are upgraded. The corresponding objective represents the lower bound (LB) for that parameter combination. 
Once all of the LBs have been computed, we find the objective value (Obj) for each budget level and report the difference as extra transportation time (ETT). The  result tables for each level of the covered population ($p$ value) are located in \ref{Ap-results}, Tables \ref{p24}-\ref{p9}. The values for Obj and ETT are reported in minutes.

Table \ref{extra} presents the average objective value for each combination of budget and $p$ levels. The first row and the first column specify the specific settings. In the first row, for each $p$ value, the average LB is also reported. In each column, we report the average value of Obj, ETT, and percentage of ETT reduction (improvement) compared to the minimum budget level, i.e., 5\% of $\hat{B}$, (column ``Imp"). Recall that $\hat{B}$ is the total budget required to upgrade all vulnerable road segments.

\vspace{-10pt}
\begin{table}[!htb]
	\tiny
	\renewcommand{\arraystretch}{1}
	\centering
	\caption{Reduction in ETT}
	\begin{adjustbox}{max width=\textwidth}	
		\begin{tabular}{lrrrrrrrrrrrrrrr}
			\hline\hline
			& \multicolumn{3}{c}{(p=24, LB=104,786)} &       & \multicolumn{3}{c}{(p=18, LB=112,236)} &       & \multicolumn{3}{c}{(p=14, LB=120,209)} &       & \multicolumn{3}{c}{(p=9, LB=127,216)} \\
			\cmidrule{2-4}\cmidrule{6-8}\cmidrule{10-12}\cmidrule{14-16}    Budget (of $\hat{B}$) & \multicolumn{1}{c}{Obj} & \multicolumn{1}{c}{ETT} & Imp   &       & \multicolumn{1}{c}{Obj} & \multicolumn{1}{c}{ETT} & Imp   &       & \multicolumn{1}{c}{Obj} & \multicolumn{1}{c}{ETT} & Imp   &       & \multicolumn{1}{c}{Obj} & \multicolumn{1}{c}{ETT} & Imp \\
			5.0\% &          106,029  &          1,244  & -     &       &          114,031  &          1,795  & -     &       &          124,487  &          4,278  & -     &       &          133,485  &          6,270  & - \\
			7.5\% &          105,157  &             372  & 70\%  &       &          112,660  &             425  & 76\%  &       &          120,813  &             604  & 81\%  &       &          127,974  &             758  & 84\% \\
			10.0\% &          104,807  &                21  & 98\%  &       &          112,275  &                39  & 98\%  &       &          120,337  &             128  & 96\%  &       &          127,412  &             197  & 97\% \\
			15.0\% &          104,786  &                  0  & 100\% &       &          112,236  &                  0  & 100\% &       &          120,209  &                  0  & 100\% &       &          127,216  &                  0  & 100\% \\
			\hline\hline
		\end{tabular}%
	\end{adjustbox}
	\label{extra}%
\end{table}%

\vspace{-10pt}
From Table \ref{extra},  an increase in the budget will result in reductions in the total transportation time in the network until it reaches the LB. Most of the reduction in ETT  occurs
when the budget level is increased from 5\% to 7.5\% and 10\% of $\hat{B}$, respectively. 
Finally, the remaining ETT is eliminated by setting $B=15\%$ of $\hat{B}$. This shows that a mitigation budget equal to 15\% of $\hat{B}$ will be enough to improve the performance of Coralville's transportation network to that of its fully functional status.  
This does not mean that all of the vulnerable roads have been upgraded, but with those that have been upgraded, the travel times of the residents from population centers to HCFs become the same as they were before the flood. 

\textbf{Insight:} A mitigation budget equal to 15\% of $\hat{B}$ will be enough to improve the performance of Coralville's transportation network up to its fully functional status for getting to the HCFs.

 The average number and length of the upgraded roads at incremental budget levels are reported in Table \ref{numlen}. We see that, on \emph{average}, an increase (decrease) in the amount of available budget  leads  to an increase (decrease) in the number of upgraded roads. But this is not always the case.  To illustrate this, we select Instance 44 from Table \ref{p18} and run additional experiments by reducing the budget from 15\% to 10\% in 10 steps (i.e., with a step-size=0.5\%). More roads are upgraded when the budget is at 11\% or 11.5\% levels than 12\% in Figure \ref{BneM}.  In these situations, the RNFMP solution usually changes to replace a long road segment consisting of multiple lanes with several short roads consisting of one or two lanes. As a result, even though a larger \textit{number} of roads get upgraded, the total length of the newly added roads can be shorter than the length of the replaced ones. This proves that increasing the budget and upgrading more \textit{miles} of road segments do not always yield a larger \textit{number} of upgraded roads. The numbers inside the brackets report the length of the upgraded roads (in miles) for each budget level.

\textbf{Insight:} Increasing the budget and upgrading more \textit{miles} of road segments do not always yield a larger \textit{number} of upgraded roads.

\vspace{-10pt}
\begin{minipage}[t]{0.55\textwidth}
	\centering
		\renewcommand{\arraystretch}{1.5}
	\captionof{table}{The average number and length of the upgraded roads}
\begin{tabular}{lcc}
	\hline\hline
	\makecell{Budget \\ (of $\hat{B}$)} & \makecell{Upgraded segments \\ (\# roads)} &  \makecell{Upgraded segments \\ (miles)} \\
	\hline
	5\%   & 38 & 11.24 \\
	7.5\% & 57 & 18.10 \\
	10\%  & 68 & 23.20 \\
	15\%  & 84 & 26.84 \\
	\hline\hline
\end{tabular}%
\label{numlen}%
\end{minipage}
\begin{minipage}[t]{0.4\textwidth}
	\centering
	\captionof{figure}{The number of upgraded roads as the budget decreases from 15\% to 10\% of $\hat{B}$}
	\includegraphics[scale=0.35]{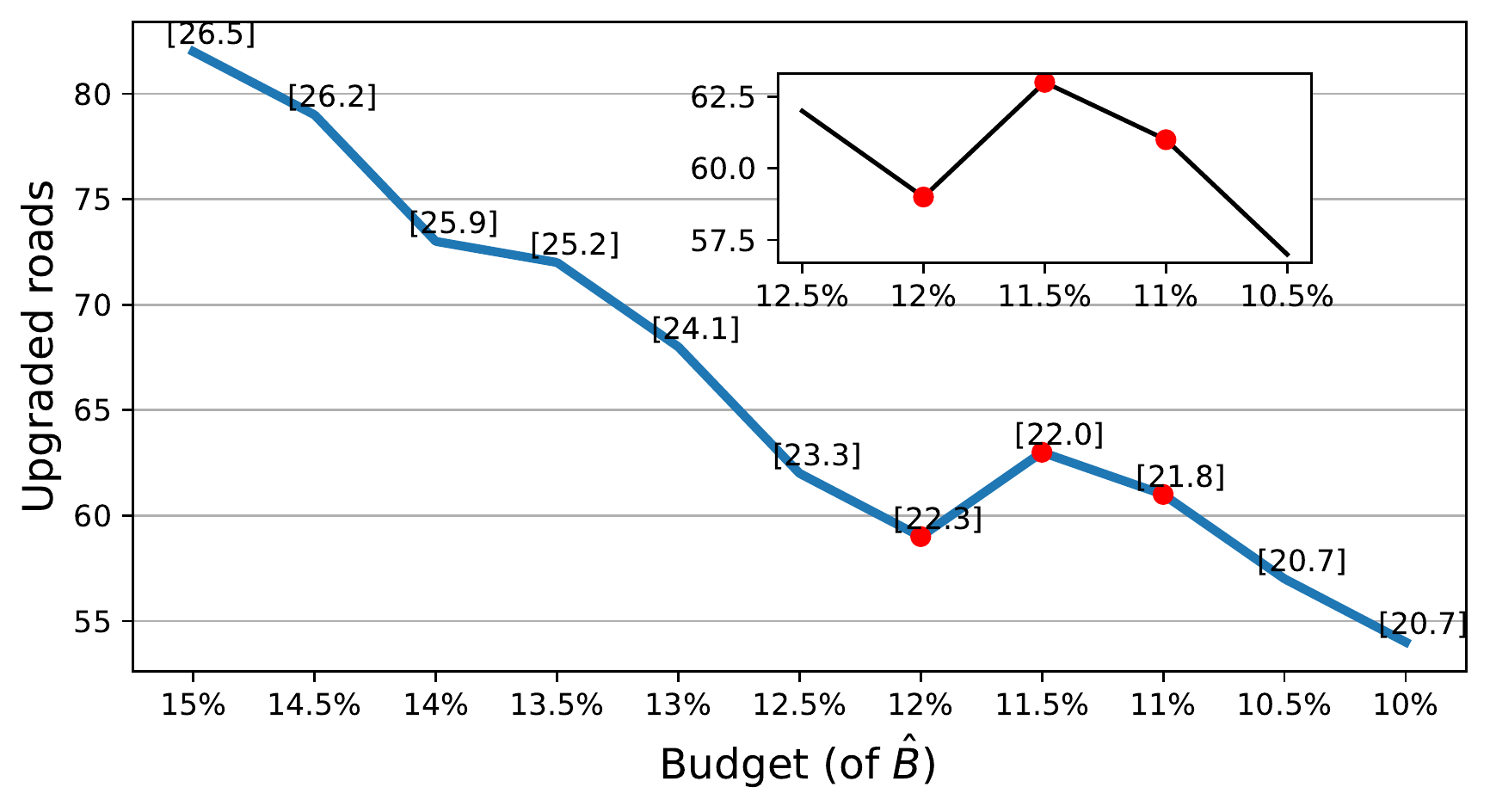}
\label{BneM}
\end{minipage}

\subsection{High frequency roads for upgrading}\label{high}
In total, 158 out of the existing 268 vulnerable roads are upgraded across our 128 instances. The optimal set of vulnerable roads for upgrading can vary quite a bit depending on the particular parameter setting, especially the value of $p$. This makes sense because instances with smaller values of $p$ will have more origins than instances with larger values of $p$.

\begin{figure}[h!]
	\centering
		\caption{The location of the upgraded roads with 100\% frequency}
	\includegraphics[scale=0.25]{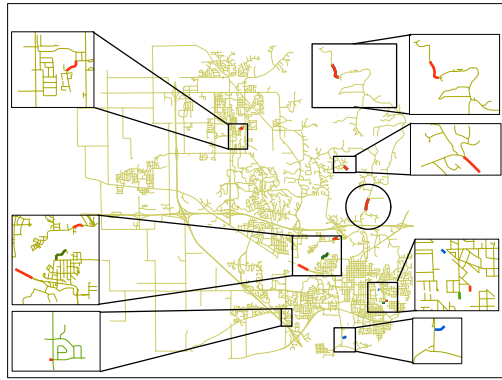}
	\label{NewFreqAnls}
\end{figure}

We divide our 128 instances into four subsets of instances based on the four different $p$ values ($4\times32$). For each subset, we compute the number of times each road has been upgraded over the 32 instances and compute the associated percentage. The roads upgraded in 100\% of the instances in each subset are reported in Table \ref{p_priority}. The first seven columns report the features of these roads. The last four columns show the subsets of instances where the roads are selected based on the $p$ values. A checkmark indicates that the corresponding road is upgraded in 100\% of the instances in that particular subset. 
Each of the roads selected in 100\% of the instances is required to keep all of the origins connected to the network.  Thus, in Table \ref{p_priority}, with smaller values of $p$ and a larger number of origins, the number of upgraded roads either remains the same or increases. When $p=24$ and $p=18$, the set of upgraded roads (checkmarks) are the same, indicating that eight roads are upgraded in all 128 instances. 
Figure \ref{NewFreqAnls} presents the location of these roads, colored in red. The roads upgraded in instances with $p=14$ and $p=9$ (only the last two columns) are in green. Finally, the ones upgraded only in instances with $p=9$  (the last column) are colored in blue. 

\textbf{Insight:}
The set of upgraded roads for Coralville always includes eight vulnerable roads that are needed to create connectivity between the network's OD pairs.

\vspace{-10pt}
\begin{table}[htbp]
	\tiny
	\centering
	\caption{The roads with 100\% frequency for each $p$ value}
			\renewcommand{\arraystretch}{0.9}
	\begin{adjustbox}{max width=\textwidth}		
		\begin{tabular}{lrccccrlcccc}
			\hline\hline
			Type  & \multicolumn{1}{l}{Length (ft)} & \multicolumn{1}{l}{Speed Limit (mi)} & \multicolumn{1}{l}{\# lanes} & Oneway & Contains bridge &       & Name  & $p=24$  & $p=18$  & $p=14$  & $p=9$ \\
			\hline
			Tertiary & 881   & 35    & 2     & no    & no    &       & Sugar Bottom Road Northeast & $\checkmark$     & $\checkmark$     & $\checkmark$     & $\checkmark$ \\
			Tertiary & 688   & 45    & 2     & no    & no    &       & Prairie du Chien Road Northeast & $\checkmark$     & $\checkmark$     & $\checkmark$     & $\checkmark$ \\
			Residential & 298   & 30    & 2     & no    & no    &       & Stewart Road & $\checkmark$     & $\checkmark$     & $\checkmark$     & $\checkmark$ \\
			Residential & 297   & 25    & 2     & no    & no    &       & Foster Road & $\checkmark$     & $\checkmark$     & $\checkmark$     & $\checkmark$ \\
			Residential & 289   & 25    & 2     & no    & no    &       & Manor Drive &       &       & $\checkmark$     & $\checkmark$ \\
			Residential & 244   & 25    & 2     & no    & no    &       & Commercial Drive &       &       &       & $\checkmark$ \\
			Residential & 242   & 25    & 2     & no    & yes   &       & Rachael Street & $\checkmark$     & $\checkmark$     & $\checkmark$     & $\checkmark$ \\
			Residential & 171   & 25    & 2     & no    & no    &       & Normandy Drive &       &       & $\checkmark$     & $\checkmark$ \\
			Residential & 161   & 25    & 2     & yes   & no    &       & Idyllwild Drive & $\checkmark$     & $\checkmark$     & $\checkmark$     & $\checkmark$ \\
			Residential & 116   & 25    & 2     & no    & no    &       & South 7th Avenue &       &       & $\checkmark$     & $\checkmark$ \\
			Residential & 98    & 25    & 2     & no    & no    &       & Manor Drive &       &       & $\checkmark$     & $\checkmark$ \\
			Residential & 95    & 25    & 2     & no    & no    &       & 6th Avenue & $\checkmark$     & $\checkmark$     & $\checkmark$     & $\checkmark$ \\
			Secondary & 84    & 25    & 2     & no    & no    &       & Muscatine Avenue &       &       &       & $\checkmark$ \\
			Residential & 49    & 35    & 4     & no    & no    &       & Rushmore Drive & $\checkmark$     & $\checkmark$     & $\checkmark$     & $\checkmark$ \\
			\hline\hline
		\end{tabular}%
	\end{adjustbox}
	\label{p_priority}%
\end{table}%
\vspace{-10pt}
To better present the role of these roads in increasing the network's connectivity, the vulnerable roads that either were not upgraded or were upgraded in less than $100\%$ of the instances are removed from Figure \ref{NewFreqAnls}. That is why some roads look disconnected (see Figure \ref{flooded} to see the full network representing the location of all vulnerable roads). Figure \ref{NewFreqAnls} helps demonstrate how these roads, except for the circled one, called Prairie du Chien Road Northeast, are necessary to prevent any  nodes from being disconnected from the network. For example, the blue road located on the southeast part of the network should be upgraded to guarantee the connectivity of its endpoint with the rest of the network. However, the red road located on the northeast part of the network provides connectivity for more nodes than just its endpoint.

\begin{figure}[h!]
		\centering
	\caption{Two alternative paths connecting the residents located at the northeast side of the network to MC-H}
\begin{minipage}[c]{0.43\textwidth}
\includegraphics[width=\linewidth]{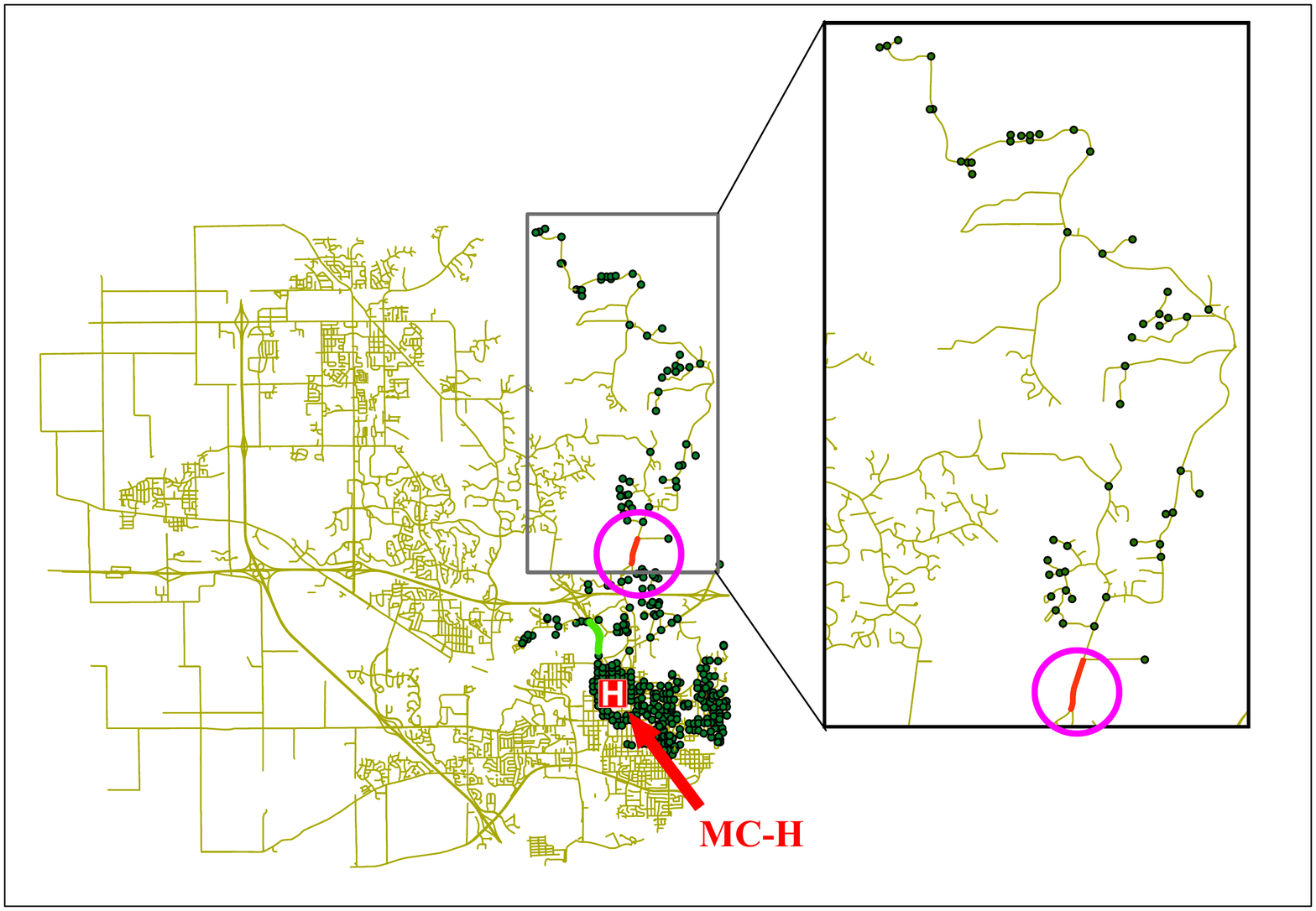}
\end{minipage}
\begin{minipage}[c]{0.4\textwidth}
	\includegraphics[width=0.97\linewidth]{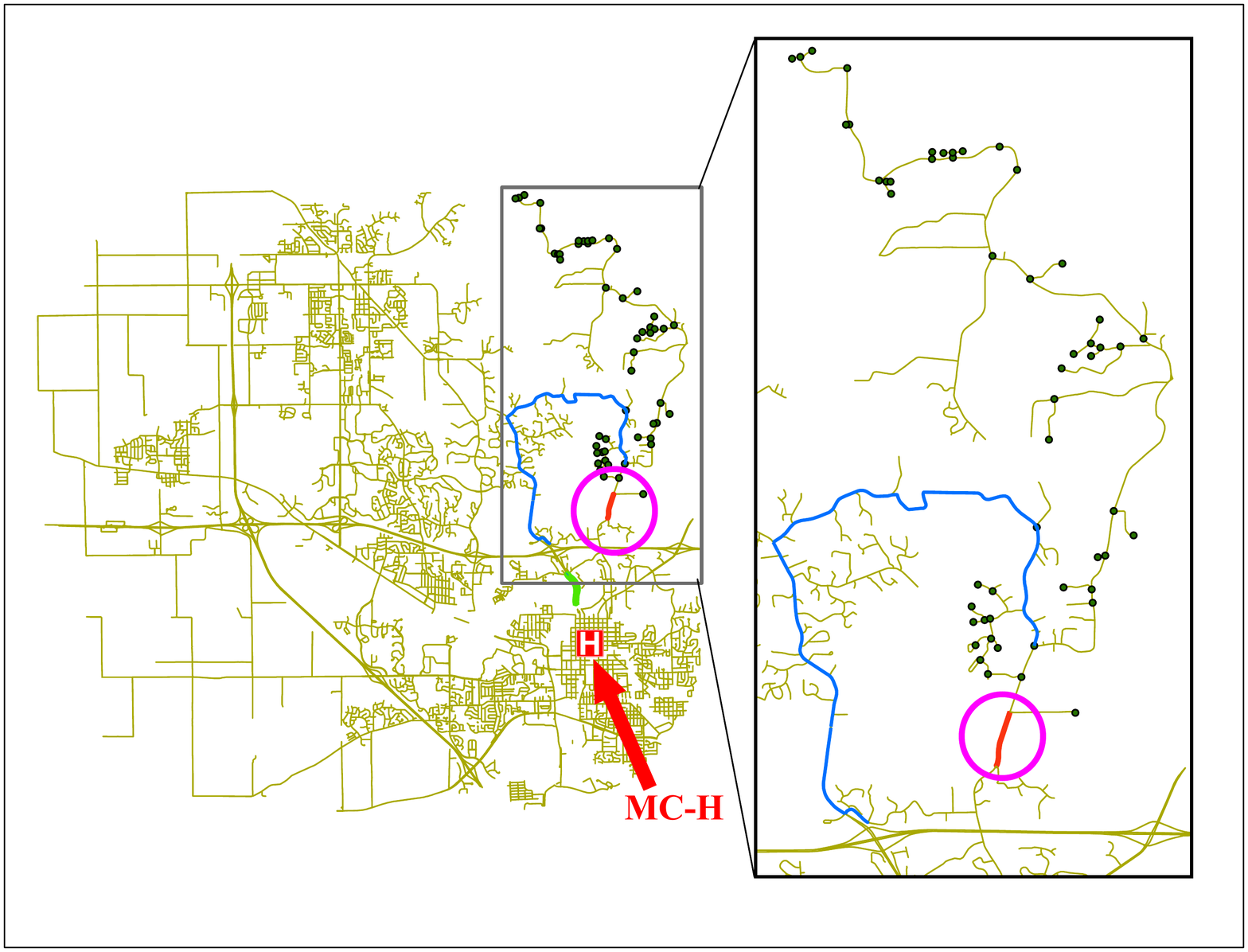}
\end{minipage}
	\label{inst13}
\end{figure}

To realize why the circled road in Figure \ref{NewFreqAnls} has also been selected to be upgraded in all instances, we choose Instance 13 from Table \ref{p24} and study the results in more detail. Figure \ref{inst13} (left panel) shows the population centers assigned to the MC-H facility (black nodes). The residents located on the northeast part of the network traverse this specific road to reach MC-H in the optimal solution.
The right panel shows a potential alternative path colored in blue where the corresponding residents have to choose to access the facility if the red road is not upgraded. 
In this example, it is easy to see how upgrading this road can significantly improve the travel time of a large number of population centers compared to using the blue path. More specifically, the alternative path will increase the total weighted travel time of the corresponding residents by 6,856.76. Similarly, in the next paragraph, we will calculate this extra weighted travel time (EWTT) added to the network's shortest paths by not upgrading every single vulnerable road.

To better understand the roads that are selected with high frequency, but are not necessary to satisfy the network's connectivity, we will focus  on Instance 13 (with $p=24$ and $|\mathcal{N}_d|=3$) as an example. To get a better understanding of what makes an arc valuable, we find the weighted shortest paths between every existing OD pair in the fully functional network (weighted by the population at the origin). Next, we iteratively choose a road from the set of vulnerable roads (not including those necessary to retain connectivity). If the selected road is located on any of the shortest paths, we remove it from the network and compute the length of the alternative shortest paths. By comparing the length of the alternative shortest paths with the original ones, we can characterize the value of a vulnerable road. We will refer to this difference as EWTT.  

Table \ref{20ETT} shows the top 20 roads among all vulnerable ones ranked based on their corresponding EWTT values. 
The first column shows the rank of each road based on EWTT values represented in the third column. The name of each road is reported in the second column. The third column specifies the upgrading cost of each road. Finally, in the last four columns, a checkmark ($\checkmark$) means that the corresponding road was upgraded when the budget was equal to $15\%$, $10\%$, $7.5\%$ or $5\%$ of the $\hat{B}$.

The first three roads (rows) are three different sections of North Dubuque Street, the green road located above MC-H in Figure \ref{NewFreqAnls}. Note that they are three \textit{connected segments} of the same road on the network. Therefore, upgrading only one of them cannot create a connected path. Thus, they should be either upgraded altogether or none of the three. That's why either all of them are upgraded (when $B=7.5\%,10\%,$ and $15\%$) or none of them (when $B=5\%$). The OpenStreetMap ID (OSMID) of each road segment, available in \cite{DOT2019}, is provided in front of each road's name inside the parenthesis. In transportation networks, roads must be split into segments when attributes change, e.g., speed limit, number of lanes, width, surface condition, etc. Any changes in these attributes can impact the performance of the network. Thus, OSMIDs are used to label each segment of the roads. In RNFMP, a change in the speed limit or the number of lanes can impact the road's upgrading cost and EWTT values. In case the change does not correspond to these two attributes, using the properties of Technique 7, these segments will be considered as one road in the model. However, they will be considered as separate segments when we report the results.

The fourth row in Table \ref{20ETT} shows Prairie du Chien Road Northeast, which is the same red (circled) road shown in Figure \ref{NewFreqAnls}. It can be seen that this road has the highest EWTT value among the ones that are upgraded in all four budget values.

\vspace{-12pt}
\begin{table}[htbp]
	\tiny
	\centering
	\caption{The top 20 roads that removing them adds the largest EWTT to the network}
			\renewcommand{\arraystretch}{1}
	\begin{adjustbox}{max width=\textwidth}			
		\begin{tabular}{lllrccrrcccc}
			\hline\hline
			Rank  & Name  & Type  & Length (ft) & \# lanes & Contains bridge & EWTT  & Upgrading cost (\$) & B=15\% & B=10\% & B=7.5\% & B=5\% \\
			\hline
			1     & North Dubuque Street (43919429)   & Secondary &          711.7  & 2     & no    & 11,688 & 28,388 & \checkmark & \checkmark & \checkmark &  \\
			2     & North Dubuque Street (15945464)  & Secondary &            48.4  & 2     & no    & 11,688 & 1,929 & \checkmark & \checkmark & \checkmark &  \\
			3     & North Dubuque Street (43919429)  & Secondary &          203.7  & 2     & no    & 10,544 & 8,124 & \checkmark & \checkmark & \checkmark &  \\
			4     & Prairie du Chien Road Northeast (15930807) & Tertiary &          688.3  & 2     & yes   & 6,856 & 27,455 & \checkmark & \checkmark & \checkmark & \checkmark \\
			5     & 2nd Street (653146696)  & Primary &            14.9  & 1     & no    & 6,082 & 298   & \checkmark & \checkmark & \checkmark & \checkmark \\
			6     & 1st Avenue (43912551)  & Secondary &          981.2  & 2     & no    & 4,748 & 39,139 & \checkmark & \checkmark & \checkmark & \checkmark \\
			7     & 1st Avenue (15937581)  & Secondary &          258.9  & 2     & no    & 4,748 & 5,164 & \checkmark & \checkmark & \checkmark & \checkmark \\
			8     & 2nd Street (621495156)  & Primary &          488.5  & 1     & no    & 3,475 & 9,741 & \checkmark & \checkmark & \checkmark & \checkmark \\
			9     & 2nd Street (15932245) & Primary &          285.9  & 1     & no    & 3,120 & 5,703 & \checkmark & \checkmark & \checkmark & \checkmark \\
			10    & Oakdale Boulevard (15941202) & Residential &          700.1  & 2     & yes   & 2,756 & 27,926 & \checkmark &       &       &  \\
			11    & Oakdale Boulevard (595674874) & Residential &          174.8  & 2     & no    & 2,756 & 6,971 & \checkmark &       &       &  \\
			12    & 2nd Street (15935233) & Primary &      1,028.8  & 1     & no    & 2,107 & 20,518 & \checkmark & \checkmark & \checkmark &  \\
			13    & Grand Army of the Republic Highway (184474115) & Primary &      2,233.3  & 3     & no    & 1,241 & 133,625 &       & \checkmark &       &  \\
			14    & North Dubuque Street (15945442) & Secondary &          210.0  & 2     & no    & 1,129 & 8,375 & \checkmark & \checkmark & \checkmark &  \\
			15    & Mormon Trek Boulevard, 1st Avenue (15938073,15938190)& Secondary &          631.5  & 2     & no    & 901   & 25,190 & \checkmark &       & \checkmark &  \\
			16    & 1st Avenue (15938073) & Secondary &          243.1  & 3     & no    & 671   & 14,543 & \checkmark &       & \checkmark &  \\
			17    & Old Hwy 218 South (15933691)& Primary &          443.5  & 2     & no    & 395   & 17,689 & \checkmark & \checkmark & \checkmark &  \\
			18    & Muscatine Avenue (15946428) & Secondary &          174.6  & 2     & yes   & 351   & 6,964 & \checkmark & \checkmark & \checkmark & \checkmark \\
			19    & Muscatine Avenue (15946428) & Secondary &            73.0  & 2     & no    & 338   & 2,913 & \checkmark & \checkmark & \checkmark & \checkmark \\
			20    & West Zeller Street (15946665) & Residential &          122.8  & 2     & no    & 276   & 4,898 & \checkmark & \checkmark &       &  \\
			\hline\hline
		\end{tabular}%
	\end{adjustbox}
	\label{20ETT}%
\end{table}%
\vspace{-12pt}

The fifth column shows that when there is enough budget (i.e., $B=15\%$), all of these roads, except for road 13, are selected to be upgraded. To find out why road 13 is not upgraded, we can look at the upgrading cost of these roads (fourth column). It can be seen that road 13 is significantly more expensive than others for upgrading. Additionally, the last four columns show that as we decrease the budget, the number of upgraded roads selected from the top 20 ones decreases. The reason is that either there is insufficient budget to upgrade those roads or the RNFMP should replace them with cheaper roads ranked below 20.   

\textbf{Insight:}
In large mitigation budgets for Coralville, the roads 
on either long OD paths, or the paths from origins with the most population
are added to the set of upgraded roads.
Removing them from Coralville's network increases the length of ODs' weighted shortest paths the most.

			\subsection{Sensitivity analysis of the number of available HCFs}  
			In this section, we divide our instances into two subsets of instances based on the number of available HCFs and summarize the insights from using seven HCFs as compared with three facilities. For each subset, we present the average objective value (Obj) and the corresponding average ETT achieved at each budget level. Remember that the ETT values report the gap between Obj and LB values where LB values are calculated after solving RNFMP over the fully functional network, i.e., by ignoring the budget constraints. The Imp sub-columns in Table \ref{tab:HCFs} show the improvements in the average Obj and ETT, achieved by utilizing all seven HCFs compared with the three hospitals, equal 28\% and 39\%, respectively. Additionally, we can observe that with seven HCFs and utilizing only 5\% of $\hat{B}$, we can further reduce the total weighted travel times 
			compared to upgrading the networks with a budget equal to 15\% of the $\hat{B}$ but utilizing only three HCFs.
		Comparing the last four sub-columns provides another evidence that upgrading more roads is not necessarily equal to upgrading longer lengths of roads.

		\textbf{Insight:} In a low-budget mitigation plan,
		utilizing all seven HCFs of Coralville
		can further improve the residents' travel times
		compared to a high-budget plan in which only three HCFs are utilized. 
		
		\begin{table}[htbp]
			\tiny
			\centering
			\caption{Summary of the results based on the number of HCFs}
					\renewcommand{\arraystretch}{1}
			\begin{adjustbox}{max width=\textwidth}			
				\begin{tabular}{cccccrrrcccccc}
					\hline\hline
					& \multicolumn{3}{c}{Obj} &       & \multicolumn{3}{c}{ETT} &       & \multicolumn{2}{c}{Upgraded segments (\# roads)} &       & \multicolumn{2}{c}{Upgraded segments (miles)} \\
					\cmidrule{2-4}\cmidrule{6-8}\cmidrule{10-11}\cmidrule{13-14}    \multicolumn{1}{l}{Budget (of $\hat{B}$)} & 3 HCFs & 7 HCFs & Imp   &       & 3 HCFs & 7 HCFs & Imp   &       & 3 HCFs & 7 HCFs &       & 3 HCFs & 7 HCFs \\
					\multicolumn{1}{l}{5\%} & 139,000 & 99,925 & 28\%  &       & 4,314 & 2,389 & 45\%  &       & 36    & 40    &       & 11.33 & 11.30 \\
					\multicolumn{1}{l}{7.50\%} & 135,281 & 98,021 & 28\%  &       & 595   & 485   & 18\%  &       & 54    & 60    &       & 17.60 & 18.08 \\
					\multicolumn{1}{l}{10\%} & 134,269 & 97,635 & 27\%  &       & 99    & 97    & 2\%   &       & 67    & 68    &       & 22.32 & 22.27 \\
					\multicolumn{1}{l}{15\%} & 134,687 & 97,537 & 28\%  &       & 0     & 0     & 0\%   &       & 86    & 82    &       & 27.82 & 26.04 \\
					\hline
					\multicolumn{3}{c}{Average:} & 28\%  &       &       &       & 39\%  &       & 61    & 63    &       & 19.77 & 19.42 \\
					\hline\hline
				\end{tabular}%
			\end{adjustbox}		
			\label{tab:HCFs}%
		\end{table}%

		We can see that the largest ETT values and the corresponding improvements belong to the instances where $B=5\%$ of $\hat{B}$.
	
\begin{figure}[h!]
			\centering
						\caption{The gap distribution between the objective and target values}
			\includegraphics[scale=0.25]{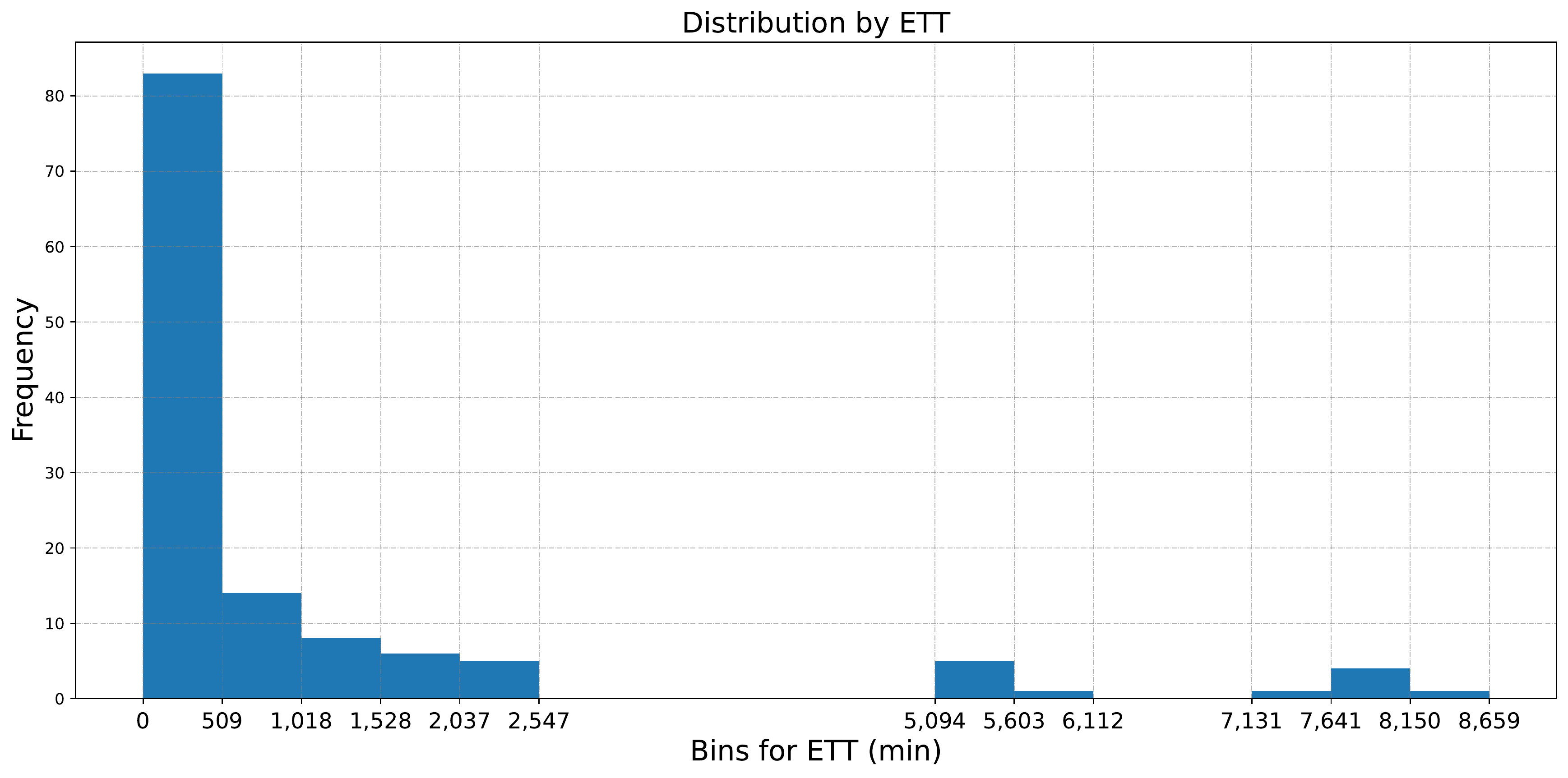}
			\label{HCF}
		\end{figure}

		The histogram depicted in Figure \ref{HCF} shows that while 116 ETT values out of the total 128 instances lie between 0 and 2,547, there are 12 points located within the range of 5,094 to 8,660. These 12 points represent the worst-case scenarios for our case study, having the largest gap between the objective value and the LB.  Additionally, Table \ref{worsttab} reports the detailed results for these cases, ranked based on the value of the gap between their objective and the  associated LB (i.e., ETT). They elucidate that the budget is at its lowest level in all worst-case scenarios, i.e., 5\% of $\hat{B}$. Next, we can see that the number of the HCFs has the most significant impact on the improvements, as the number of HCFs equals three in the top four worst-case scenarios and in eight out of the 12 worst-case scenarios. From the above, we can see that the benefits of employing smaller HCFs, such as clinics (in addition to the hospitals), for providing emergency services in the event of a flood are crucial in improving the network's performance. However, the value of the available budget has the highest impact on minimizing the residents' total ETT.

		\textbf{Insight:} To minimize the residents' total ETT in Coralville, the value of the available budget should be the main concern of the federal and local authorities, followed by the number of available HCFs.

	\vspace{-15pt}
		\begin{table}[h!]
			\tiny
			\centering
			\caption{The 12 worst-case scenarios with the largest gaps between the objective and LB values (ETT) }
					\renewcommand{\arraystretch}{1}
			\begin{tabular}{cccrcccccc}
				\hline\hline
				Rank  & Instance &       & p     & HCFs  & Capacity & $\alpha$ & B     & Obj   & ETT\\
				\hline
				1     & 109   &       & 9     & 3     & Different & 30\%  & 5\%   &          160,248  &          8,659  \\
				2     & 105   &       & 9     & 3     & Different & 15\%  & 5\%   &          161,450  &          8,088  \\
				3     & 101   &       & 9     & 3     & Identical & 30\%  & 5\%   &          150,349  &          7,943  \\
				4     & 97    &       & 9     & 3     & Identical & 15\%  & 5\%   &          150,547  &          7,784  \\
				5     & 121   &       & 9     & 7     & Different & 15\%  & 5\%   &          141,377  &          7,780  \\
				6     & 125   &       & 9     & 7     & Different & 30\%  & 5\%   &          136,844  &          7,619  \\
				7     & 77    &       & 14    & 3     & Different & 30\%  & 5\%   &          149,050  &          5,642  \\
				8     & 89    &       & 14    & 7     & Different & 15\%  & 5\%   &          131,923  &          5,556  \\
				9     & 69    &       & 14    & 3     & Identical & 30\%  & 5\%   &          140,048  &          5,334  \\
				10    & 73    &       & 14    & 3     & Different & 15\%  & 5\%   &          150,323  &          5,238  \\
				11    & 93    &       & 14    & 7     & Different & 30\  & 5\%   &          127,338  &          5,121  \\
				12    & 65    &       & 14    & 3     & Identical & 15\%  & 5\%   &          140,164  &          5,116  \\
				\hline\hline
			\end{tabular}%
			\label{worsttab}%
		\end{table}%

\subsection{Sensitivity analysis of HCF capacity}\label{cap_anls}
In Table \ref{alpha}, the first two rows 
		summarize the results corresponding to the type of capacity allocation to HCFs, i.e., the ``Identical" capacity allocation versus the ``Different" allocation policy, which is based on the number of available beds. We observe that within the identical policy, the average objective and gap values have been reduced by 21.0\% and 40.3\% compared to the ``Different" policy, respectively. In  Figure \ref{flooded}, we can see that the three hospitals with more beds than other HCFs
		are located very close to each other. Therefore, we can conclude that the network can better serve the residents, and its performance is higher when the total capacity is allocated more evenly across the city, whereas currently, the majority of the capacity is in one location.

		The impact of different $\alpha$ values is summarized in the third and fourth rows.  The results of Table \ref{alpha} show that the capacity allocation policy has a bigger influence on the network's performance than $\alpha$ values. Figure \ref{cap_alpha} further illuminates that the gap between different policies of capacity allocation is more noticeable than different $\alpha$ values. Note that when the setting changes from (Different, 30\%) to (Identical, 15\%), the capacity of each HCF decreases. Nevertheless, the network's performance improves by almost 21.2\% due to the change in the type of capacity allocation.

		\textbf{Insight:} 
		With the same total capacity among Coralville's HCFs, we can have a much better objective value if the capacity is more evenly distributed among the network. For Coralville, it can be helpful to have a large hospital (i.e., more capacity) on the north side of the network.

\begin{minipage}[t]{0.55\textwidth}
\centering
\renewcommand{\arraystretch}{1}
\captionof{table}{Summary of the results for different levels of $\alpha$\\ and types of capacity allocation}
\begin{tabular}{llcrcc}
\hline\hline
\multicolumn{2}{l}{Setting} & Obj   & ETT   & \# roads & miles \\
\hline
Capacity: & Different &                  130,964  &               1,361  &                                                         58  &                                                      19  \\
& Identical &                  103,252  &                   633  &                                                         66  &                                                      20  \\
&       &       &       &       &  \\
$\alpha$: & 15\%  &                  117,935  &                   825  &                                                         62  &                                                      20  \\
& 30\%  &                  115,981  &                   868  &                                                         62  &                                                      20  \\
\hline\hline
\end{tabular}%
\label{alpha}%
\hfill
\end{minipage}
\begin{minipage}[t]{0.45\textwidth}
\centering
\captionof{figure}{Impact of HCFs' capacity on the objective}
\includegraphics[width=0.8\textwidth]{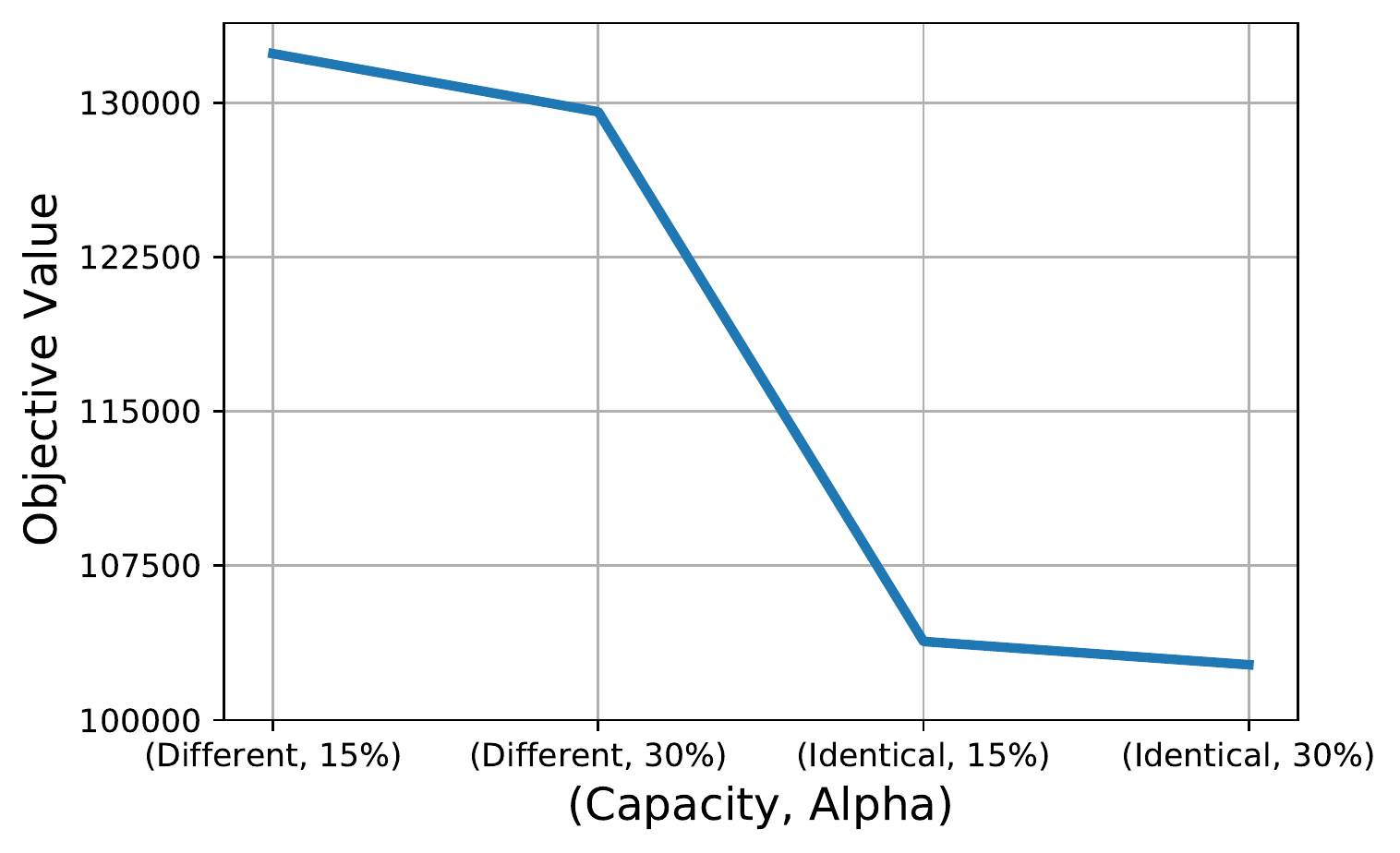}
\label{cap_alpha}
\end{minipage}

\vspace{-5pt}
\section{Conclusions and Future Work}\label{con}
In this study, we have developed a network design problem that minimizes the transportation times between population centers and HCFs by making selected upgrades on the road segments. We have modeled a transportation network that overlays on a physical road network, identified the complexity of the problem, and presented  improvements that make the problem tractable. We have demonstrated the value of this problem in a case study involving the city of Coralville and analyzed the impact each parameter has on the results.

In this study, we had access to a unique level of data and information from the road networks in the State of Iowa. While high-quality regulatory maps, like those used in this study, require significant data and computational resources, data-driven approaches \citep{li2022accounting,hu2021real} provide a high-quality and affordable alternative for flood map generation for resource-limited communities. Computationally, there is potential to improve the runtimes of the model, particularly for larger networks consisting of more nodes and arcs. Additionally, there is potential value in studying the version of the problem that minimizes the longest transportation time on the network. This involves considering the difference between the minsum and minmax approaches. 

There are also some practical issues that can be handled in future research. Incorporating the road's capacity or traffic could be an important part of the model to expand. This would require utilizing the travel time's delay function, which necessitates developing a non-linear model.
Also, 
the developed RNFMP model is based on the assumption that the locations of the HCFs are given. However, the network's performance may be quite limited by the current locations. Thus, higher accessibility gains could be provided by efficiently locating temporary facilities. We studied the role of the RNFMP in increasing accessibility to HCFs. However, it seems appropriate that future network investment models need to be developed to consider not only impacts to access to HCFs but also access to local (super)markets, gas stations, police departments, and fire departments.

\ACKNOWLEDGMENT{The authors would like to thank the Mid-America Transportation Center (MATC) and the U.S. Department of Transportation (DOT) for sponsoring this research (Grant \# 69A3551747107), and Iowa Flood Center for providing flood maps. 
}

\bibliographystyle{pomsref}

 \let\oldbibliography\thebibliography
 \renewcommand{\thebibliography}[1]{%
    \oldbibliography{#1}%
    \baselineskip10pt 
    \setlength{\itemsep}{8pt}
 }
\bibliography{ref1}

\ECSwitch %

\ECHead{E-Companion}

\section{Proof of Propositions} \label{proofs}
\subsection{Proposition \ref{yVI}}
	\begin{repeatproposition}
	If there exists a node $i\in \mathcal{N}_o$ that is not incident to any non-vulnerable outgoing arc (i.e., $\nexists j \in \mathcal{N}: (i,j) \in \mathcal{A}_n$), and $\delta_{out}(i)\ge 1$, then:
	\begin{align}		 
		&\sum_{j:(i,j) \in \mathcal{A}_v} y_{ij} \ge 1 \label{VI}
	\end{align}
\end{repeatproposition}
\proof{Proof of Proposition \ref{yVI}.} Since all of the outgoing arcs incident to $i$ are vulnerable, it is necessary to upgrade at least one of those outgoing arcs to be able to send out the flow from $i$ and have a feasible solution (Figure \ref{tech}l). In case $\delta_{out}(i)=1$, the corresponding $y_{ij}$ variable can be eliminated from RNFMP and set road $(i,j)$ as an upgraded road. Subsequently, the budget value $B$ should be updated to $B-c_{ij}$, and the corresponding weighted travel time $t_{ij}w^k$ should be added to the value of the objective function in the pre-processing step. 
\Halmos  
\endproof

\subsection{Proposition \ref{prop2}}
\begin{repeatproposition}
	For every origin $k\in \mathcal{N}_o$, which is connected to all of the existing destinations in the non-upgraded network  $(i.e., \mathbb{G}\setminus \mathcal{A}_v)$, variable $x_{ij}^k$ can be eliminated from RNFMP if
	\begin{align}		 
		&SP_f(k,i)+t_{ij}> \max_{g \in \mathcal{N}_d}\{SP_n(k,g)\} & \forall (i,j)\in \mathcal{A} \label{p2}
	\end{align}
	where $SP_n(u,v)$ and $SP_f(u,v)$ are the lengths of the shortest paths from node $u\in \mathcal{N}_o$ to node $v\in \mathcal{N}$ in the non-upgraded $(\mathbb{G}\setminus \mathcal{A}_v)$ and fully upgraded $(\mathbb{G})$, respectively.
\end{repeatproposition}
\proof{Proof of Proposition \ref{prop2}.}
The length of the shortest path between every pair of nodes $u$ and $v$ in the non-upgraded network will either remain the same or get shorter after upgrading the network. That is, $SP_f(u,v) \le SP_n(u,v)$. The shortest paths in a non-upgraded network can only contain non-vulnerable arcs. Thus, in the optimal solution, the only case in which a particular shortest path specified in the $\mathbb{G}\setminus \mathcal{A}_v$ network may get replaced by an upgraded (but costly) shortest path in the \textit{upgraded} network ($\mathbb{G}$) occurs when the costly upgraded path provides a shorter travel time and a better objective value. Accordingly, $\max_{g \in \mathcal{N}_d}\{SP_n(k,g)\}$, which is the length of the shortest path from $k$ to its furthest destination in $\mathbb{G}\setminus \mathcal{A}_v$ provides an upper bound on the travel time of the residents associated with origin $k$ to any destination. Therefore, it is not optimal for the residents associated with origin $k$ to visit arc $(i,j)\in\mathcal{A}$ if the corresponding travel time in the fully upgraded network ($\mathbb{G}$) is larger than the upper bound, i.e., $\max_{g \in \mathcal{N}_d}\{SP_n(k,g)\}$.
\Halmos  
\endproof

\subsection{Proposition \ref{delta}}
\begin{repeatproposition}
	If there exists a disconnected component $\Delta$ obtained from $\mathbb{G}\setminus\{v\}$, where $v\in \mathcal{N}_{p}$ and $\nexists u\in \mathcal{N}_{\Delta} : u \in \mathcal{N}_d$, then $x_{ij}^k$ variables can be eliminated for every
	$k \notin \mathcal{N}_{\Delta}$ and $(i,j)\in\mathcal{A}_{\Delta}$.
\end{repeatproposition}
\proof{Proof of Proposition \ref{delta}.}
Since there exists no destination node in component $\Delta$, the residents associated with origin $k \notin \mathcal{N}_{\Delta}$ who enter the component need to leave it again through node $v$.
This creates a circuit that, which is not optimal and the decision variables corresponding to such paths (i.e., $x_{ij}^k:$ $\forall k \in \mathcal{N}_o\setminus \mathcal{N}_{\Delta}$, $\forall (i,j)\in \mathcal{A}_{\Delta}$) can be eliminated from the search space. 
\Halmos  
\endproof

\section{Impact of the improvements for Coralville} \label{impact-coral}
	\begin{table}[h!]
	\renewcommand{\arraystretch}{1.4}
	\centering
	\caption{Impact of the improvements in instances of network $\mathbb{G}_C$}
	\begin{adjustbox}{max width=\textwidth}	
		\begin{tabular}{rlrrrrrrrrrrrrrr}
			\hline\hline
			\multicolumn{1}{l}{Instance} & \multicolumn{1}{c}{Budget} & \multicolumn{1}{c}{$p$} & \multicolumn{1}{c}{$\mathbb{G}$} &       & \multicolumn{2}{c}{$\mathbb{G}^\prime$} &       & \multicolumn{2}{c}{$\mathbb{G}^\prime_{IS}$} &       & \multicolumn{2}{c}{$\mathbb{G}^\prime_{IV}$} &       & \multicolumn{2}{c}{$\mathbb{G}^\prime_{IS+VI}$} \\
			\cmidrule{6-7}\cmidrule{9-10}\cmidrule{12-13}\cmidrule{15-16}          &       &       & \multicolumn{1}{c}{RT} &       & \multicolumn{1}{c}{RT} & \multicolumn{1}{c}{Imp} &       & \multicolumn{1}{c}{RT} & \multicolumn{1}{c}{Imp} &       & \multicolumn{1}{c}{RT} & \multicolumn{1}{c}{Imp} &       & \multicolumn{1}{c}{RT} & \multicolumn{1}{c}{Imp} \\
			\hline
			\multicolumn{1}{l}{1} & \multicolumn{1}{c}{15\%} & \multicolumn{1}{c}{24} & 1,020.72 &       & 448.92 & 56\%  &       & 523.19 & 49\%  &       & 511.07 & 50\%  &       & 512.53 & 50\% \\
			\multicolumn{1}{l}{2} &       & \multicolumn{1}{c}{18} & 1,491.10 &       & 644.52 & 57\%  &       & 750.37 & 50\%  &       & 821.69 & 45\%  &       & 785.15 & 47\% \\
			\multicolumn{1}{l}{3} &       & \multicolumn{1}{c}{14} & 1,931.31 &       & 833.19 & 57\%  &       & 990.00 & 49\%  &       & 997.09 & 48\%  &       & 975.97 & 49\% \\
			\multicolumn{1}{l}{4} &       & \multicolumn{1}{c}{9} & 2,678.43 &       & 1,174.88 & 56\%  &       & 1,488.56 & 44\%  &       & 1,376.54 & 49\%  &       & 1,545.43 & 42\% \\
			\hline
			\multicolumn{1}{l}{5} & \multicolumn{1}{c}{10\%} & \multicolumn{1}{c}{24} & 1,003.90 &       & 610.42 & 39\%  &       & 592.67 & 41\%  &       & 519.08 & 48\%  &       & 558.35 & 44\% \\
			\multicolumn{1}{l}{6} &       & \multicolumn{1}{c}{18} & 1,429.06 &       & 918.57 & 36\%  &       & 894.47 & 37\%  &       & 764.68 & 46\%  &       & 782.04 & 45\% \\
			\multicolumn{1}{l}{7} &       & \multicolumn{1}{c}{14} & 1,969.96 &       & 1,211.56 & 38\%  &       & 1,191.09 & 40\%  &       & 1,015.44 & 48\%  &       & 1,085.50 & 45\% \\
			\multicolumn{1}{l}{8} &       & \multicolumn{1}{c}{9} & 3,037.65 &       & 1,770.17 & 42\%  &       & 1,900.38 & 37\%  &       & 1,663.30 & 45\%  &       & 1,763.27 & 42\% \\
			\hline
			\multicolumn{1}{l}{9} & \multicolumn{1}{c}{7.5\%} & \multicolumn{1}{c}{24} & 1,258.83 &       & 556.28 & 56\%  &       & 760.73 & 40\%  &       & 652.62 & 48\%  &       & 674.23 & 46\% \\
			\multicolumn{1}{l}{10} &       & \multicolumn{1}{c}{18} & TL(0.01\%) &       & 8,120.89 & 25\%  &       & 6,751.08 & 37\%  &       & 7,084.19 & 34\%  &       & 4,521.48 & 58\% \\
			\multicolumn{1}{l}{11} &       & \multicolumn{1}{c}{14} & 5,602.01 &       & 4,707.45 & 16\%  &       & 4,234.74 & 24\%  &       & 3,762.54 & 33\%  &       & 4,620.15 & 18\% \\
			\multicolumn{1}{l}{12} &       & \multicolumn{1}{c}{9} & TL(0.03\%) &       & 3,910.46 & 64\%  &       & 4,400.74 & 59\%  &       & 5,092.05 & 53\%  &       & 4,790.74 & 56\% \\
			\hline
			\multicolumn{1}{l}{13} & \multicolumn{1}{c}{5\%} & \multicolumn{1}{c}{24} & 2,702.19 &       & 1,474.24 & 45\%  &       & 1,347.63 & 50\%  &       & 1,786.75 & 34\%  &       & 1,294.71 & 52\% \\
			\multicolumn{1}{l}{14} &       & \multicolumn{1}{c}{18} & 1,690.73 &       & 1,285.68 & 24\%  &       & 1,235.58 & 27\%  &       & 1,319.70 & 22\%  &       & 1,315.48 & 22\% \\
			\multicolumn{1}{l}{15} &       & \multicolumn{1}{c}{14} & TL(0.22\%) &       & TL(0.20\%) & 0\%   &       & TL(0.02\%) & 0\%   &       & TL(0.04\%) & 0\%   &       & TL(0.02\%) & 0\% \\
			\multicolumn{1}{l}{16} &       & \multicolumn{1}{c}{9} & TL(0.01\%) &       & 7,445.32 & 31\%  &       & 3,619.67 & 66\%  &       & 2,637.89 & 76\%  &       & 4,222.18 & 61\% \\
			\hline
			& Average: &       & 4,313.49 &       & 2,869.54 & 40\%  &       & 2,592.56 & 41\%  &       & 2,550.29 & 42\%  &       & 2,515.45 & 42\% \\
		\end{tabular}%
	\end{adjustbox}
	\label{imp-C}%
\end{table}%

\newpage
\section{Detailed results for the case study} \label{Ap-results}
\begin{table}[h!]
	\centering
	\caption{Results of the case study for $p=24$}
		\renewcommand{\arraystretch}{1.3}
	\begin{adjustbox}{max width=\textwidth}	
		\begin{tabular}{lclrrrrcc}
			\hline\hline
			Instance & \# HCFs & Capacity & \multicolumn{1}{c}{$\alpha$} & \multicolumn{1}{c}{B} & \multicolumn{1}{c}{ Obj } & \multicolumn{1}{c}{Ext} & \# upgraded roads & Length of upgrades (mi) \\
			\hline\
			1     & 3     & Identical & 15\%  & 5.0\% &    118,777  &          1,444  & 32    & 10.96 \\
			2     & 3     & Identical & 15\%  & 7.5\% &    117,436  &             103  & 57    & 18.75 \\
			3     & 3     & Identical & 15\%  & 10.0\% &    117,346  &                13  & 75    & 24.22 \\
			4     & 3     & Identical & 15\%  & 15.0\% &    117,333  &                  0  & 76    & 26.40 \\
			5     & 3     & Identical & 30\%  & 5.0\% &    118,745  &          1,806  & 32    & 10.96 \\
			6     & 3     & Identical & 30\%  & 7.5\% &    117,042  &             103  & 51    & 17.80 \\
			7     & 3     & Identical & 30\%  & 10.0\% &    116,939  &                  0  & 75    & 24.22 \\
			8     & 3     & Identical & 30\%  & 15.0\% &    116,939  &                  0  & 75    & 24.22 \\
			9     & 3     & Different & 15\%  & 5.0\% &    128,066  &          1,575  & 32    & 11.02 \\
			10    & 3     & Different & 15\%  & 7.5\% &    127,180  &             688  & 51    & 18.03 \\
			11    & 3     & Different & 15\%  & 10.0\% &    126,527  &                36  & 57    & 20.81 \\
			12    & 3     & Different & 15\%  & 15.0\% &    126,493  &                  2  & 75    & 24.59 \\
			13    & 3     & Different & 30\%  & 5.0\% &    126,584  &          1,586  & 35    & 11.14 \\
			14    & 3     & Different & 30\%  & 7.5\% &    125,667  &             669  & 52    & 18.03 \\
			15    & 3     & Different & 30\%  & 10.0\% &    125,038  &                40  & 55    & 20.90 \\
			16    & 3     & Different & 30\%  & 15.0\% &    124,998  &                  0  & 76    & 25.08 \\
			17    & 7     & Identical & 15\%  & 5.0\% &      69,257  &             191  & 36    & 11.66 \\
			18    & 7     & Identical & 15\%  & 7.5\% &      69,074  &                  7  & 62    & 18.55 \\
			19    & 7     & Identical & 15\%  & 10.0\% &      69,066  &                  0  & 66    & 21.22 \\
			20    & 7     & Identical & 15\%  & 15.0\% &      69,066  &                  0  & 66    & 21.22 \\
			21    & 7     & Identical & 30\%  & 5.0\% &      67,445  &             189  & 36    & 11.58 \\
			22    & 7     & Identical & 30\%  & 7.5\% &      67,263  &                  7  & 62    & 18.55 \\
			23    & 7     & Identical & 30\%  & 10.0\% &      67,257  &                  0  & 66    & 21.22 \\
			24    & 7     & Identical & 30\%  & 15.0\% &      67,257  &                  1  & 65    & 21.22 \\
			25    & 7     & Different & 15\%  & 5.0\% &    111,524  &          1,583  & 33    & 11.14 \\
			26    & 7     & Different & 15\%  & 7.5\% &    110,643  &             702  & 51    & 18.29 \\
			27    & 7     & Different & 15\%  & 10.0\% &    109,985  &                44  & 55    & 20.78 \\
			28    & 7     & Different & 15\%  & 15.0\% &    109,944  &                  2  & 76    & 25.05 \\
			29    & 7     & Different & 30\%  & 5.0\% &    107,835  &          1,578  & 33    & 11.14 \\
			30    & 7     & Different & 30\%  & 7.5\% &    106,950  &             693  & 51    & 18.09 \\
			31    & 7     & Different & 30\%  & 10.0\% &    106,293  &                36  & 56    & 20.75 \\
			32    & 7     & Different & 30\%  & 15.0\% &    106,257  &                  0  & 75    & 24.59 \\
			\hline\hline
		\end{tabular}%
	\end{adjustbox}
	\label{p24}%
\end{table}%

\begin{table}[H]
	\centering
	\caption{Results of the case study for $p=18$}
		\renewcommand{\arraystretch}{1.3}
	\begin{adjustbox}{max width=\textwidth}	
		\begin{tabular}{lclrrrrcc}
			\hline\hline
			Instance & \# HCFs & Capacity & \multicolumn{1}{c}{$\alpha$} & \multicolumn{1}{c}{B} & \multicolumn{1}{c}{ Obj } & \multicolumn{1}{c}{ETT} & \# upgraded roads & Length of upgrades (mi) \\
			\hline\
			33    & 3     & Identical & 15\%  & 5.0\% &    127,898  &          2,021  & 33    & 11.05 \\
			34    & 3     & Identical & 15\%  & 7.5\% &    126,023  &             145  & 54    & 18.12 \\
			35    & 3     & Identical & 15\%  & 10.0\% &    125,889  &                11  & 81    & 25.31 \\
			36    & 3     & Identical & 15\%  & 15.0\% &    125,878  &                  0  & 84    & 28.45 \\
			37    & 3     & Identical & 30\%  & 5.0\% &    127,898  &          2,340  & 33    & 11.05 \\
			38    & 3     & Identical & 30\%  & 7.5\% &    125,758  &             199  & 52    & 17.57 \\
			39    & 3     & Identical & 30\%  & 10.0\% &    125,559  &                  0  & 81    & 25.31 \\
			40    & 3     & Identical & 30\%  & 15.0\% &    125,558  &                  0  & 83    & 26.27 \\
			41    & 3     & Different & 15\%  & 5.0\% &    137,676  &          2,181  & 30    & 10.92 \\
			42    & 3     & Different & 15\%  & 7.5\% &    136,246  &             750  & 51    & 18.04 \\
			43    & 3     & Different & 15\%  & 10.0\% &    135,560  &                65  & 54    & 20.71 \\
			44    & 3     & Different & 15\%  & 15.0\% &    135,495  &                  0  & 82    & 26.52 \\
			45    & 3     & Different & 30\%  & 5.0\% &    136,185  &          2,272  & 30    & 10.92 \\
			46    & 3     & Different & 30\%  & 7.5\% &    134,647  &             734  & 51    & 18.04 \\
			47    & 3     & Different & 30\%  & 10.0\% &    133,987  &                74  & 55    & 20.68 \\
			48    & 3     & Different & 30\%  & 15.0\% &    133,913  &                  0  & 83    & 27.00 \\
			49    & 7     & Identical & 15\%  & 5.0\% &      73,926  &             302  & 41    & 11.25 \\
			50    & 7     & Identical & 15\%  & 7.5\% &      73,640  &                16  & 69    & 18.49 \\
			51    & 7     & Identical & 15\%  & 10.0\% &      73,624  &                  0  & 76    & 23.64 \\
			52    & 7     & Identical & 15\%  & 15.0\% &      73,624  &                  0  & 76    & 23.64 \\
			53    & 7     & Identical & 30\%  & 5.0\% &      72,027  &             314  & 41    & 11.25 \\
			54    & 7     & Identical & 30\%  & 7.5\% &      71,731  &                18  & 64    & 18.44 \\
			55    & 7     & Identical & 30\%  & 10.0\% &      71,713  &                  0  & 77    & 23.79 \\
			56    & 7     & Identical & 30\%  & 15.0\% &      71,713  &                  0  & 76    & 23.79 \\
			57    & 7     & Different & 15\%  & 5.0\% &    120,316  &          2,492  & 31    & 10.92 \\
			58    & 7     & Different & 15\%  & 7.5\% &    118,590  &             766  & 52    & 18.07 \\
			59    & 7     & Different & 15\%  & 10.0\% &    117,906  &                82  & 56    & 20.60 \\
			60    & 7     & Different & 15\%  & 15.0\% &    117,826  &                  1  & 84    & 27.07 \\
			61    & 7     & Different & 30\%  & 5.0\% &    116,321  &          2,442  & 31    & 10.92 \\
			62    & 7     & Different & 30\%  & 7.5\% &    114,650  &             770  & 52    & 18.07 \\
			63    & 7     & Different & 30\%  & 10.0\% &    113,959  &                80  & 55    & 20.58 \\
			64    & 7     & Different & 30\%  & 15.0\% &    113,885  &                  6  & 82    & 26.52 \\
			\hline\hline
		\end{tabular}%
	\end{adjustbox}
	\label{p18}%
\end{table}%

\newpage
\begin{table}[ht]
	\renewcommand{\arraystretch}{1}
	\centering
	\caption{Results of the case study for $p=14$}
		\renewcommand{\arraystretch}{1.3}
	\begin{adjustbox}{max width=\textwidth}	
		\begin{tabular}{lclrrrrcc}
			\hline\hline
			Instance & \# HCFs & Capacity & \multicolumn{1}{c}{$\alpha$} & \multicolumn{1}{c}{B} & \multicolumn{1}{c}{ Obj } & \multicolumn{1}{c}{ETT} & \# upgraded roads & Length of upgrades (mi) \\
			\hline\
			65    & 3     & Identical & 15\%  & 5.0\% &    140,164  &          5,116  & 38    & 11.76 \\
			66    & 3     & Identical & 15\%  & 7.5\% &    135,313  &             265  & 53    & 17.61 \\
			67    & 3     & Identical & 15\%  & 10.0\% &    135,087  &                39  & 80    & 24.16 \\
			68    & 3     & Identical & 15\%  & 15.0\% &    135,048  &                  0  & 94    & 30.27 \\
			69    & 3     & Identical & 30\%  & 5.0\% &    140,048  &          5,334  & 38    & 11.76 \\
			70    & 3     & Identical & 30\%  & 7.5\% &    135,112  &             399  & 53    & 17.61 \\
			71    & 3     & Identical & 30\%  & 10.0\% &    134,733  &                20  & 78    & 24.13 \\
			72    & 3     & Identical & 30\%  & 15.0\% &    134,714  &                  0  & 90    & 28.09 \\
			73    & 3     & Different & 15\%  & 5.0\% &    150,323  &          5,238  & 36    & 11.66 \\
			74    & 3     & Different & 15\%  & 7.5\% &    146,105  &          1,020  & 51    & 17.03 \\
			75    & 3     & Different & 15\%  & 10.0\% &    145,292  &             207  & 52    & 19.84 \\
			76    & 3     & Different & 15\%  & 15.0\% &    145,086  &                  0  & 89    & 28.33 \\
			77    & 3     & Different & 30\%  & 5.0\% &    149,050  &          5,642  & 36    & 11.69 \\
			78    & 3     & Different & 30\%  & 7.5\% &    144,395  &             986  & 53    & 17.26 \\
			79    & 3     & Different & 30\%  & 10.0\% &    143,596  &             187  & 55    & 20.17 \\
			80    & 3     & Different & 30\%  & 15.0\% &    143,409  &                  0  & 90    & 28.81 \\
			81    & 7     & Identical & 15\%  & 5.0\% &      79,157  &             733  & 45    & 11.34 \\
			82    & 7     & Identical & 15\%  & 7.5\% &      78,450  &                26  & 70    & 18.34 \\
			83    & 7     & Identical & 15\%  & 10.0\% &      78,424  &                  0  & 82    & 25.22 \\
			84    & 7     & Identical & 15\%  & 15.0\% &      78,424  &                  0  & 82    & 25.22 \\
			85    & 7     & Identical & 30\%  & 5.0\% &      76,918  &             510  & 46    & 11.41 \\
			86    & 7     & Identical & 30\%  & 7.5\% &      76,439  &                31  & 70    & 18.34 \\
			87    & 7     & Identical & 30\%  & 10.0\% &      76,409  &                  0  & 81    & 25.22 \\
			88    & 7     & Identical & 30\%  & 15.0\% &      76,409  &                  0  & 82    & 25.22 \\
			89    & 7     & Different & 15\%  & 5.0\% &    132,894  &          6,527  & 36    & 11.69 \\
			90    & 7     & Different & 15\%  & 7.5\% &    127,423  &          1,056  & 52    & 16.93 \\
			91    & 7     & Different & 15\%  & 10.0\% &    126,687  &             321  & 61    & 20.68 \\
			92    & 7     & Different & 15\%  & 15.0\% &    126,367  &                  0  & 89    & 28.39 \\
			93    & 7     & Different & 30\%  & 5.0\% &    127,338  &          5,121  & 36    & 11.66 \\
			94    & 7     & Different & 30\%  & 7.5\% &    123,269  &          1,053  & 50    & 16.89 \\
			95    & 7     & Different & 30\%  & 10.0\% &    122,468  &             251  & 55    & 19.96 \\
			96    & 7     & Different & 30\%  & 15.0\% &    122,219  &                  2  & 88    & 28.09 \\
			\hline\hline
		\end{tabular}%
	\end{adjustbox}
	\label{p14}%
\end{table}%

\begin{table}[ht]
	\renewcommand{\arraystretch}{1}
	\centering
	\caption{Results of the case study for $p=9$}
		\renewcommand{\arraystretch}{1.3}
	\begin{adjustbox}{max width=\textwidth}	
		\begin{tabular}{lclrrrrcc}
			\hline\hline
			Instance & \# HCFs & Capacity & \multicolumn{1}{c}{$\alpha$} & \multicolumn{1}{c}{B} & \multicolumn{1}{c}{ Obj } & \multicolumn{1}{c}{ETT} & \# upgraded roads & Length of upgrades (mi) \\
			\hline\
			97    & 3     & Identical & 15\%  & 5.0\% &    150,547  &          7,784  & 44    & 11.67 \\
			98    & 3     & Identical & 15\%  & 7.5\% &    143,227  &             464  & 61    & 17.30 \\
			99    & 3     & Identical & 15\%  & 10.0\% &    142,805  &                42  & 80    & 24.02 \\
			100   & 3     & Identical & 15\%  & 15.0\% &    142,763  &                  0  & 97    & 31.48 \\
			101   & 3     & Identical & 30\%  & 5.0\% &    150,349  &          7,943  & 44    & 11.67 \\
			102   & 3     & Identical & 30\%  & 7.5\% &    142,980  &             575  & 62    & 17.31 \\
			103   & 3     & Identical & 30\%  & 10.0\% &    142,441  &                35  & 82    & 24.06 \\
			104   & 3     & Identical & 30\%  & 15.0\% &    142,405  &                  0  & 96    & 29.30 \\
			105   & 3     & Different & 15\%  & 5.0\% &    161,450  &          8,088  & 45    & 11.49 \\
			106   & 3     & Different & 15\%  & 7.5\% &    154,572  &          1,210  & 56    & 15.01 \\
			107   & 3     & Different & 15\%  & 10.0\% &    153,734  &             372  & 61    & 20.60 \\
			108   & 3     & Different & 15\%  & 15.0\% &    153,362  &                  0  & 96    & 29.95 \\
			109   & 3     & Different & 30\%  & 5.0\% &    160,248  &          8,659  & 45    & 11.49 \\
			110   & 3     & Different & 30\%  & 7.5\% &    152,792  &          1,203  & 61    & 18.03 \\
			111   & 3     & Different & 30\%  & 10.0\% &    151,943  &             354  & 61    & 19.72 \\
			112   & 3     & Different & 30\%  & 15.0\% &    151,593  &                  4  & 97    & 30.35 \\
			113   & 7     & Identical & 15\%  & 5.0\% &      84,455  &             998  & 52    & 11.11 \\
			114   & 7     & Identical & 15\%  & 7.5\% &      83,521  &                64  & 69    & 18.24 \\
			115   & 7     & Identical & 15\%  & 10.0\% &      83,464  &                  6  & 88    & 25.77 \\
			116   & 7     & Identical & 15\%  & 15.0\% &      83,458  &                  0  & 93    & 29.19 \\
			117   & 7     & Identical & 30\%  & 5.0\% &      82,138  &             813  & 51    & 10.98 \\
			118   & 7     & Identical & 30\%  & 7.5\% &      81,389  &                64  & 69    & 18.25 \\
			119   & 7     & Identical & 30\%  & 10.0\% &      81,331  &                  6  & 88    & 25.77 \\
			120   & 7     & Identical & 30\%  & 15.0\% &      81,325  &                  0  & 92    & 29.04 \\
			121   & 7     & Different & 15\%  & 5.0\% &    141,853  &          8,256  & 44    & 11.41 \\
			122   & 7     & Different & 15\%  & 7.5\% &    134,844  &          1,247  & 58    & 17.84 \\
			123   & 7     & Different & 15\%  & 10.0\% &    133,974  &             378  & 60    & 20.52 \\
			124   & 7     & Different & 15\%  & 15.0\% &    133,600  &                  3  & 89    & 29.14 \\
			125   & 7     & Different & 30\%  & 5.0\% &    136,844  &          7,619  & 44    & 11.38 \\
			126   & 7     & Different & 30\%  & 7.5\% &    130,465  &          1,240  & 58    & 17.84 \\
			127   & 7     & Different & 30\%  & 10.0\% &    129,606  &             381  & 61    & 20.62 \\
			128   & 7     & Different & 30\%  & 15.0\% &    129,225  &                  0  & 93    & 29.21 \\
			\hline\hline
		\end{tabular}%
	\end{adjustbox}
	\label{p9}%
\end{table}%

\newpage

\section{Sensitivity analysis of the covered population}\label{pop_anals}
The total length and number of times that each type of road has been upgraded over each $p$ value is presented in the first four rows of Table \ref{rtype_p}. The fifth row shows that for covering larger populations, more (miles of) residential roads need to be upgraded. However, the number of upgraded segments among the other types of roads and the number of upgraded roads containing bridges (last row) are almost insensitive to the percentage of the covered population. Table \ref{p_priority} in Section \ref{high} lists the roads with the 100\% frequency for each subset of instances with a specific $p$ value. In total, the number of upgraded roads with 100\% frequency where $p=$ 24, 18, 14, and 9 is equal to 9, 9, 12, and 14 roads, respectively. This
indicates that for covering larger populations more upgraded roads are needed to connect the OD pairs.

\begin{table}[htbp]
	\centering
	\caption{The number and length of the upgraded roads over different $p$ values}
					\renewcommand{\arraystretch}{1}
	\begin{tabular}{llrrrrrrrrrrr}
		\hline\hline
		&       & \multicolumn{2}{c}{p=24} &       & \multicolumn{2}{c}{p=18} &       & \multicolumn{2}{c}{p=14} &       & \multicolumn{2}{c}{p=9} \\
		\cmidrule{3-4}\cmidrule{6-7}\cmidrule{9-10}\cmidrule{12-13}    Type  &       & roads & miles &       & roads & miles &       & roads & miles &       & roads & miles \\
		\hline\
		Primary &       & 13    &      73,755  &       & 13    &      72,258  &       & 13    &      69,413  &       & 13    &      64,270  \\
		Secondary &       & 17    &      83,812  &       & 18    &      81,599  &       & 18    &      74,123  &       & 18    &      66,276  \\
		Tertiary &       & 12    &      68,042  &       & 15    &      69,933  &       & 13    &      68,354  &       & 18    &      73,050  \\
		Residential &       & 66    &    144,768  &       & 81    &    162,844  &       & 91    &    191,258  &       & 104   &    222,229  \\
		Total &       & 108   &    370,377  &       & 127   &    386,634  &       & 135   &    403,148  &       & 153   &    425,825  \\
		\multicolumn{2}{l}{Roads with bridge} & 16    &      70,855  &       & 16    &      68,798  &       & 16    &      64,971  &       & 18    &      68,681  \\
		\hline\hline
	\end{tabular}%
	\label{rtype_p}%
\end{table}%

\end{document}